\documentclass[ 12pt]{article}

\relax
\usepackage{xspace,mathrsfs}

\def\tableofcontents{\section*{Contents\@mkboth{CONTENTS}{CONTENTS}}
	\@starttoc{toc}}
\def\l@part#1#2{\addpenalty{\@secpenalty}
	\addvspace{2.25em plus 1\p@} \begingroup
	\@tempdima 3em \parindent \z@ \rightskip \@pnumwidth \parfillskip
	-\@pnumwidth
	{\large \textbf \leavevmode #1\hfil \hbox to\@pnumwidth{\hss #2}}\par
	\nobreak \endgroup}
\def\l@section#1#2{\addpenalty{\@secpenalty} \addvspace{0em}
	\@tempdima 1.5em \begingroup
	\parindent \z@ \rightskip \@pnumwidth
	\parfillskip -\@pnumwidth
	\textbf \leavevmode \advance\leftskip\@tempdima \hskip -\leftskip #1\nobreak\hfil
	\nobreak\hbox to\@pnumwidth{\hss #2}\par
	\endgroup}

\relax
\usepackage{xspace,mathrsfs}

\usepackage[nottoc]{tocbibind}
\usepackage{mathrsfs}
\usepackage{comment}
\usepackage[polutonikogreek,english]{babel}
\usepackage{calligra}
\usepackage{amsmath, amssymb}
\usepackage{amsthm}
\usepackage{xcolor}
\usepackage{soul}
\usepackage{scalerel}
\usepackage[only,llbracket,rrbracket,llparenthesis,rrparenthesis]{stmaryrd}
\usepackage{accsupp} 
\usepackage{tikz}
\usepackage{tikz-cd}
\newsavebox{\pullback}
\sbox\pullback{%
	\begin{tikzpicture}%
		\draw (0,0) -- (1ex,0ex);%
		\draw (1ex,0ex) -- (1ex,1ex);%
\end{tikzpicture}}
\usepackage{multirow}
 \usetikzlibrary{decorations.pathmorphing} 

\newcommand{\Pd}[2]{#1:#2\op\to \mathsf{InfSL}}

\newcommand{\ang}[1]{\langle #1 \rangle}

\newcommand{\binprod}[3]{#1\overset{p_1}{\leftarrow} #2 \overset{p_2}{\rightarrow} #3}
\newcommand{\wbinprod}[3]{#1\overset{\mathsf{p}_1}{\leftarrow} #2 \overset{\mathsf{p}_2}{\rightarrow} #3}

\newcommand{\weakbinprodg}[4]{#2\overset{\mathsf{#1}_1}{\leftarrow} #3 \overset{\mathsf{#1}_2}{\rightarrow} #4}
\newcommand{\weakbinprodga}[5]{#2\overset{\mathsf{#1}^{#5}_1}{\leftarrow} #3 \overset{\mathsf{#1}^{#5}_2}{\rightarrow} #4}
\newcommand{\weakbinprodgp}[4]{#2\overset{\mathsf{#1}'_1}{\leftarrow} #3 \overset{\mathsf{#1}'_2}{\rightarrow} #4} 

\newcommand{\bfeq}[2]{\delta^{(\mathsf{#2}_{1},\mathsf{#2}_{2})}}
\newcommand{\bfeqa}[3]{\delta^{(\mathsf{#2}^{#3}_{1},\mathsf{#2}^{#3}_{2})}}
\newcommand{\bfeqp}[2]{\delta^{(\mathsf{#2}'_{1},\mathsf{#2}'_{2})}} 
\newcommand{\CC}{\mathscr{C}}
\newcommand{\DD}{\mathscr{D}}
\newcommand{\Pdocc}[2]{#1:#2\op\to \mathsf{InfSL}}
\newcommand{\Pdoc}{\mathsf{P}:\mathscr{C}\op\to \mathsf{InfSL}}

\newcommand{\Pstrict}{\mathsf{P}^{\mathpzc{s}}}
\newcommand{\strict}{\mathpzc{s}}
\newcommand{\CCs}{\mathscr{C}_{\mathpzc{s}}}
\newcommand{\lel}[3]{#1\le #2\le #3}
\newcommand{\fs}[1]{{#1}_{\mathpzc{s}}}
\newcommand{\M}{\mathcal{M}}

\newcommand{\lrfloor}[1]{\lfloor #1 \rceil}
\newcommand{\lrbrace}[1]{\llbrace #1 \rrbrace}

\newcommand{\std}{\textbf{Std}}
\newcommand{\ml}{\textbf{ML}}
\newcommand{\PP}{\mathsf{P}}
\newcommand{\pp}{\mathsf{p}}
\newcommand{\op}{\ensuremath{^{\operatorname{op}}}}
\newcommand{\po}{\ensuremath{_{\operatorname{po}}}}
\newcommand{\hotop}{\mathsf{Ho}(\mathbf{Top})}
\newcommand{\set}{\mathsf{Set}}
\newcommand{\Gsets}{\set^G}
\newcommand{\setsG}{\set_G}

\newcommand*{\llbrace}{%
	\BeginAccSupp{method=hex,unicode,ActualText=2983}%
	\textnormal{\usefont{OMS}{lmr}{m}{n}\char102}%
	\mathchoice{\mkern-4.05mu}{\mkern-4.05mu}{\mkern-4.3mu}{\mkern-4.8mu}%
	\textnormal{\usefont{OMS}{lmr}{m}{n}\char106}%
	\EndAccSupp{}%
}
\newcommand*{\rrbrace}{%
	\BeginAccSupp{method=hex,unicode,ActualText=2984}%
	\textnormal{\usefont{OMS}{lmr}{m}{n}\char106}%
	\mathchoice{\mkern-4.05mu}{\mkern-4.05mu}{\mkern-4.3mu}{\mkern-4.8mu}%
	\textnormal{\usefont{OMS}{lmr}{m}{n}\char103}%
	\EndAccSupp{}%
}

\DeclareFontFamily{OT1}{pzc}{}
\DeclareFontShape{OT1}{pzc}{m}{it}{<->s*[1.1]pzcmi7t}{}
\DeclareMathAlphabet{\mathpzc}{OT1}{pzc}{m}{it}

\newcounter{dummy}
\usepackage[hidelinks]{hyperref}
\usepackage{enumitem}
\usepackage{cleveref}
\crefname{enumi}{}{}
\makeatletter
\newcommand\myitem[1][]{\item[#1]\refstepcounter{dummy}\def\@currentlabel{#1}}
\makeatother
\crefname{dummy}{}{}

\newtheorem*{notation}{Notation}
\newtheorem{thm}{Theorem}[section]
\newtheorem{dfn}[thm]{Definition}
\newtheorem{prop}[thm]{Proposition}
\newtheorem{lemma}[thm]{Lemma}
\newtheorem{cor}[thm]{Corollary}
\newtheorem{rmk}[thm]{Remark}
\newtheorem{example}[thm]{Example}
\newtheorem{examples}[thm]{Examples}
\newtheorem{obs}[thm]{Observation}

\usepackage{blindtext}
\usepackage{bussproofs}

\usepackage[numbers]{natbib}
\begin{document}

	\title{BIASED ELEMENTARY DOCTRINES AND QUOTIENT COMPLETIONS}
	\author{Cipriano Junior Cioffo\thanks{The author would like to acknowledge Jacopo Emmenegger, Milly Maietti and Pino Rosolini (who suggested the name ``biased elementary doctrines") for all the significant comments and fruitful discussions on this work.}\\ \small Università di Pisa}
	\date{}
	\maketitle
	
	\begin{abstract}
In this work, we fill the gap between the elementary quotient
completion introduced by Maietti and Rosolini and the exact completion
of a category with weak finite limits, as described by Carboni and
Vitale. To achieve this, we generalize Lawvere’s elementary doctrines
to apply to categories with weak finite products, referring to these
structures as biased elementary doctrines. We present two main
constructions: the first, called strictification, produces an
elementary doctrine from a biased one, while the second is an
extension of the elementary quotient completion that generalizes the
exact completion of a category with weak finite limits, even when weak
finite products are involved.
    \end{abstract}

\section{Introduction}\label{section:Introduction}
The quotient completion is a pervasive construction in mathematics and
logic that has been extensively studied in category theory. The first
explicit description of the free exact category on a category with
finite limits was provided in \cite{carbonimagno}. Later in \cite{carboni1998regular}, the result
was generalized for categories with weak finite limits. In \cite{maietti2013quotient}, Maietti and Rosolini introduced the elementary quotient completion to
provide an abstract account of the quotient model in \cite{maietti2009minimalist}. To
achieve this, they relativized the notions of equivalence relation and
quotient for Lawvere’s elementary doctrines, i.e., suitable functors
of the form $\PP: \CC\op\to \mathsf{Pos}$ from
categories with finite products into the category of posets and
order-preserving functions \cite{lawvere1969adjointness, lawvere1970equality}. According to \cite{maietti2013quotient,maietti2013elementary,maietti2015unifying},  the elementary quotient completion has been shown to generalize
the exact completion of a category with weak finite limits, provided
that finite products exist.

A primary example illustrating the application of elementary quotient completion is the category of setoids arising from the intensional level of the \textit{Minimalist Foundation} \cite{maietti2009minimalist}, see \Cref{examples: eqc}. This example is the key example of an elementary quotient completion that is not an exact completion.

Conversely, given a group $G$, the category of $G$-sets $\Gsets$, is equivalent to the exact completion of the category of free actions $\setsG$. This serves as an example of exact completion that is not an instance of elementary quotient completion, due to the absence of finite products in $\setsG$, see \Cref{examples: w prod } and \Cref{rmk: exact completion that are not elementary quotient completion}.

In this work, we aim to extend the elementary quotient completion so that it applies to categories with weak finite products, providing a comprehensive generalization of the exact completion for categories with weak finite limits.\

To achieve this, we first define biased elementary doctrines as suitable functors of the form $\PP:\CC\op\to \mathsf{Pos}$ from categories with weak finite products. The structure of a biased elementary doctrine is similar to the usual one, but the properties are restated to account for a kind of \textit{bias} due to the weak universal property of weak finite products. {For biased elementary structures, the main novelty is the need to distinguish between the equality predicate and component-wise equality predicate, see \Cref{rmk: proof-relevant/irrelevant equality}.} As expected, biased elementary doctrines generalize Lawvere's elementary doctrines.

For these structures, we identify a class of elements (called proof-irrelevant) within the fibers of weak finite products, see \Cref{dfn:proof-irrelevant elements}, which are used to derive the two main constructions.

 The first construction, called strictification, associates to every biased elementary doctrine $\PP: \CC\op\to \mathsf{Pos}$ an elementary doctrine ${\PP^\strict}:{\CCs}\op\to \mathsf{Pos}$ on the finite product completion $\CCs$ of $\CC$. The second is a quotient completion which extends both the elementary quotient completion and the exact completion of categories with weak finite limits, even in case of weak finite products. The second construction is a quotient completion that extends both the elementary quotient completion and the exact completion of categories with weak finite limits, even in cases involving weak finite products. For this construction, we prove a universal property similar to that discussed in \cite{carboni1998regular}.

In addition, we provide a slice construction for suitable biased elementary doctrines, namely those $\PP: \CC\op\to \mathsf{Pos}$ such that $\CC$ has weak pullbacks. For every object $A\in\CC$, we define the slice doctrine $\mathsf{P}_{/A}:(\CC/A)\op\to \mathsf{Pos}$ and prove that this operation commutes with the quotient completion. {Slice doctrines are a distinctive feature of the biased framework. Indeed, the slice construction of several elementary doctrines is inherently biased, as it happens in the case of slice doctrines arising from the {intensional level} of the Minimalist Foundation.}

 {In future studies of the properties of the quotient completion, which have already begun in \cite{cioffophd}, slice doctrines will play a fundamental role concerning the local cartesian closure of the quotient category.}

In conclusion, our primary example of quotient completion falls outside both exact completion and elementary quotient completion. It is given by the quotient completion of the slice doctrines arising from the intensional level of the Minimalist Foundation, which provide a description of \textit{families of setoids} in the context of doctrines.


In section \ref{section: Preliminaries} we recall the main notions about the elementary doctrines and the relationship between the elementary quotient completion and the exact completion. In section \ref{section:biased elementary doctrines}, we define biased elementary doctrines and their basic properties. Proof-irrelevant elements are discussed in section \ref{section:proof-irrelevant elements}. In section \ref{section:strictification}, we provide the strictification of a biased elementary doctrine and in section \ref{section: quotient completion} we describe the quotient completion and its universal property. In section \ref{section: concluding remarks}, we briefly recall some properties of the exact and elementary quotient completions and explore potential further developments.


\section{Preliminaries}\label{section: Preliminaries}
In this section, we recall the main notions about elementary doctrines, mainly following \cite{maietti2013quotient,maietti2013elementary}. We also refer to \cite{lawvere1969adjointness,lawvere1970equality,jacobs1999categorical}. In particular, we recall the construction of the elementary quotient completion introduced in \cite{maietti2013quotient} and its relationship to the exact completion of a category with weak finite limits presented in \cite{carbonimagno,carboni1998regular}.

The categorical structure, which we refer to simply as \textit{doctrine}, is given by a functor of the form
\[\PP:\CC\op\to \mathsf{Pos}\]
from a category $\CC$ with finite products to the category $\mathsf{Pos}$ of  partially-ordered sets (posets) and order-preserving functions. We shall refer to the value of the functor $\PP(X)$, as the \textit{fiber} of $\PP$ on $X$.

Throughout the text, we will use the term \textit{strict} finite products to refer to finite products in order to distinguish them from \textit{weak} finite products, which will be introduced in the next section. Similarly, when we wish to emphasize that the base category of a doctrine possesses strict finite products, we will refer to it as a strict doctrine.

\begin{dfn}\label{dfn:primarydoctrine}\normalfont
	Let $\CC$ be a category with strict finite products. A \textit{primary doctrine} is a functor $\mathsf{P}:\CC\op\to \mathsf{Pos}$ which takes value in the category $\mathsf{InfSL}$ of inf-semilattices and order preserving functions which preserve finite meets, i.e.
	\begin{enumerate}
		\item\label{item:primarydoctrine-1}  for every object $X\in\CC$, $\mathsf{P}(X)$ has finite meets,
		\item\label{item:primarydoctrine-2} for every arrow $f:X\to Y$ of $\CC$, the functor $\mathsf{P}_f:\mathsf{P}(Y)\to \mathsf{P}(X)$ preserves finite meets.
	\end{enumerate}
\end{dfn}


\begin{dfn}\label{dfn:elementarydoctrine}\normalfont
	Let $\CC$ be a category with strict finite products. A primary doctrine $\Pdoc$ is called \textit{elementary} if, for every object $X\in\CC$,
	there exists an element $\delta_X\in \mathsf{P}(X\times X)$ such that:
	\begin{enumerate}
		\item\label{item:elementarydoctrine-1} For every element $\alpha\in \mathsf{P}(X)$, the assignment $$\exists_{\Delta_X}(\alpha):= \mathsf{P}_{p_1}(\alpha)\wedge\delta_X$$ 
		is the left adjoint to the functor $\mathsf{P}_{\Delta_X}:\mathsf{P}(X\times X)\to \mathsf{P}(X)$.
		\item\label{item:elementarydoctrine-2} For every object $Y\in\CC$ and arrow $e:=\langle p_1,p_2,p_2\rangle: X\times Y\to X\times Y\times Y$, the assignment 
		$$\exists_e(\alpha):=\mathsf{P}_{\langle p_1,p_2\rangle}(\alpha)\wedge\mathsf{P}_{\langle p_2,p_3\rangle}(\delta_Y)$$
		for $\alpha$ in $\mathsf{P}(X\times Y)$ is left adjoint to
		$\mathsf{P}_e:\mathsf{P}(X\times Y\times Y)\to \mathsf{P}(X\times Y)$.
	\end{enumerate}
\end{dfn}
We will refer to the element $\delta_{X}\in \mathsf{P}(X\times X)$ as the \textit{fibered equality} of $X$.

 Next, we present an equivalent definition of elementary doctrine that is discussed in \cite[Remark 2.3]{maietti2013elementary}. For a proof of the following result see \cite{emmpasqroscoalgebras} or Appendix A of \cite{cioffophd}. Before proceeding, we first recall the definition of \textit{descent data}.

	\begin{dfn}\label{dfn:descentobjects}\normalfont
	Let $\Pdoc$ be a primary doctrine. If $\beta\in \mathsf{P}(X\times X)$, then ${\mathpzc{Des}}({\beta})$ is the sub-order of elements $\alpha\in \mathsf{P}(X)$ satisfying
	\begin{equation}\label{eqdes}
		\mathsf{P}_{p_1}(\alpha)\wedge\beta\le \mathsf{P}_{p_2}(\alpha).
	\end{equation}
	\end{dfn}

\begin{prop}\label{prop:equiv-elementarydoctrine}\normalfont
	Let $\CC$ be a category with strict finite products. A primary doctrine $\Pdoc$ is \textit{elementary} if and only if for every object $X\in\CC$,
	there exists an element $\delta_X\in \mathsf{P}(X\times X)$ such that:
	\begin{enumerate}
		\item\label{item:equiv-elementarydoctrine-1} $\top_X \le \mathsf{P}_{\Delta_X}(\delta_X)$.
		\item\label{item:equiv-elementarydoctrine-2}  $\mathsf{P}(X)= {\mathpzc{Des}}({\delta_{X}})$.
		\item\label{item:equiv-elementarydoctrine-3}  $\delta_{X}\boxtimes\delta_{Y}\le\delta_{X\times Y}$, where $\delta_{X}\boxtimes\delta_{Y}:= \mathsf{P}_{\langle p_1,p_3\rangle}\delta_{X}\wedge \mathsf{P}_{\langle p_2,p_4\rangle}\delta_Y$
	\end{enumerate}
where the arrows $p_i$, for $\lel{1}{i}{4}$, are the projections of the product $X\times Y\times X\times Y$.
\qed\end{prop}

We will refer to $\delta_{X}\boxtimes\delta_{Y}$ as the \textit{component-wise equality} of $X\times Y$

\begin{dfn}\label{dfn:existential-ED}\normalfont
	An elementary doctrine $\Pdoc$ is called \textit{existential} if, for every pair of objects $X_1,X_2\in\CC$ the functors $\mathsf{P}_{{p}_i}:\mathsf{P}(X_i)\to \mathsf{P}(X_1\times X_2)$, for $i=1,2$, have left adjoints $\exists_{{p}_i}:\mathsf{P}(X_1 \times X_2)\to \mathsf{P}(X_i)$ which satisfy 
	\begin{itemize}
		\item  the \textit{Beck-Chevalley} condition:
		for the pullback diagram 
		\begin{equation}\label{diagram:pullbacks existential}
		\begin{tikzcd}
			X_1\times Y \arrow[swap]{d}{1_{X_1}\times f}\arrow[]{r}{{p}_2} & Y \arrow[]{d}{f}\\
			X_1\times X_2 \arrow[swap]{r}{{p}_2} & X_2
		\end{tikzcd}	\end{equation}
		the canonical arrow $\exists_{\mathsf{p}_2}\circ \mathsf{P}_{{1_{X_1}\times f}}(-) \le \mathsf{P}_{f}\circ \exists_{\mathsf{p}_2}(-)$ is an isomorphism. The analogous condition holds for ${p}_1$.\\
		\item the \textit{Frobenius reciprocity}: for any projection $p_i:X_1\times X_2{\to} X_i$, element $\alpha\in \mathsf{P}(X_i)$, and $\beta\in \mathsf{P}(X_1\times X_2)$, the canonical arrow $\exists_{{p}_i}(\mathsf{P}_{{p}_i}\alpha \wedge \beta)\le\alpha\wedge\exists_{{p}_i}\beta$ is an isomorphism.
	\end{itemize}
\end{dfn}

\begin{dfn}\label{dfn:universal-ED}\normalfont
	An elementary doctrine $\Pdoc$ is called \textit{universal} if for every pair of objects $X_1,X_2\in\CC$, the functors $\mathsf{P}_{{p}_i}:\mathsf{P}(X)\to \mathsf{P}(X_1\times X_2)$, for $i=1,2$, have right adjoints $\forall_{{p}_i}:\mathsf{P}(X_1\times X_2)\to \mathsf{P}(X_i)$ which satisfy
	\begin{itemize}
		\item  the \textit{Beck-Chevalley} condition:
		for the pullback diagram as in (\ref{diagram:pullbacks existential})
		the canonical arrow $\mathsf{P}_{f}\circ \forall_{{p}_2}(-)\le\forall_{{p}_2}\circ \mathsf{P}_{1_{X_1}\times f}(-)$ is an isomorphism. The analogous condition holds for $\mathsf{p}_1$.
	\end{itemize}
\end{dfn}
\begin{dfn}\label{dfn:implicational}\normalfont
	A primary doctrine $\Pdoc$ is called \textit{implicational} if for every object $X\in\CC$ and element $\alpha\in \mathsf{P}(X)$ the functor $\alpha\wedge-:\mathsf{P}(X)\to \mathsf{P}(X)$ has a right adjoint $\alpha\Rightarrow-:\mathsf{P}(X)\to \mathsf{P}(X)$.
\end{dfn}

	\begin{dfn}\label{dfn:ed-comprehensions}\normalfont
	Let $\Pdoc$ be an primary doctrine and let  $X$ be an object of $\CC$. A \textit{comprehension} of an element $\alpha\in \mathsf{P}(X)$ is an arrow $\llbrace \alpha \rrbrace: C\to X$ such that $\top_C\le \mathsf{P}_{\llbrace \alpha \rrbrace}\alpha$ and which satisfies the following universal property: for every arrow $f:Y\to X$ such that $\top_Y\le \mathsf{P}_f(\alpha)$, there exists an arrow $h:\nolinebreak Y\to\nolinebreak C$ such that the following diagram commutes
	\begin{center}
		\begin{tikzcd}
			C \arrow[dashed]{r}{h} \arrow[swap]{dr}{\llbrace \alpha \rrbrace} & Y \arrow[]{d}{f}\\
			& X.
		\end{tikzcd}
	\end{center}
	The comprehension $\lrbrace{\alpha}$ is called \textit{strict} if the induced arrow $h$ is unique. When $h$ is not unique, the comprehension $\lrbrace{\alpha}$ is called \textit{weak}. A comprehension $\llbrace \alpha \rrbrace: C\to X$ is called \textit{full} if $\alpha\le \beta$ whenever $\top_C\le \mathsf{P}_{\llbrace \alpha \rrbrace}\beta$.
\end{dfn}

\begin{dfn}\label{dfn:ed-comprehensive diagonals}\normalfont
	An elementary doctrine $\Pdoc$ has \textit{comprehensive diagonals} if for any pair of arrows $f,g: A \to X$ in $\CC$, it is  $$\top_A \leq \mathsf{P}_{\langle f, g \rangle}\delta_X\ \text{ if and only if}\ f=g.$$
\end{dfn}

\begin{examples}\label{examples: elementary doctrines}\normalfont
Most of the following examples of elementary doctrines are discussed in details in \cite{maietti2013quotient}.
	\begin{itemize}
\item\label{example:ed-strictsubobject}\normalfont
If $\CC$ is a category with finite limits, one can consider the functor
\begin{equation*}
	\Pd{\mathsf{Sub}_{\CC}}{\CC}.
\end{equation*}
which sends an object of $X\in\CC$ to the poset of subobjects of $X$. The elements of $\mathsf{Sub}_\CC(X)$ will be denoted with $\lrfloor{g}$, for a monomorphism $g:Z\to X$. The action of $\mathsf{Sub}_{\CC}$ on an arrow $f:Y\to X$ sends a subobject $\lrfloor{g}$ to the subobject given by the pullback of $g$ along $f$. The fibered equalities are given by the equivalence class of the diagonals $\Delta_X$ and the doctrine has full strict comprehensions and comprehensive diagonals. If $\CC$ is regular, $\mathsf{Sub}_{\CC}$ is existential, see \cite{HUGHES2003156}.
\item If $\CC$ is a category with weak finite limits and strict finite products,  one can consider the functor
\begin{equation*}
	\Psi_{\CC}:\CC\op\to \mathsf{InfSL}
\end{equation*} 
which sends an object $X\in\CC$ to the poset reflection  $(\CC/X)\po$ of the slice category $\CC/X$. The elements of $\Psi_\CC(X)$ are called \textit{weak subobjects} or \textit{variations} in \cite{grandis2000weak} and will be denoted with $\lrfloor{g}$, for an arrow $g:Z\to X$. In \cite{palmgren2004categorical}, they are also called \textit{pre-subobjects}. The action of $\Psi_{\CC}$ on an arrow $f:Y\to X$ sends a variation $\lrfloor{g}$ to the variation given by a weak pullback of $g$ along $f$. The fiber equality is given by the equivalence class of the diagonal $\Delta_X$. This doctrine has full weak comprehensions and it is existential, with left adjoint to $\Psi_{\CC}(f)$  given by post-composition with $f$.
\item\label{example:ed-of-Hsets} If $H$ is an inf-semilattice, one can consider the functor 
\[H:\set\op\to \mathsf{InfSL}\]
which sends a set $X$ to $H^X$ ordered pointwisely, i.e.\ for $\alpha,\beta:X\to H$ it holds $\alpha\le \beta$ if and only if $\alpha(x)\le \beta(x)$ for every $x\in X$. The action of $H$ on a function $f:Y\to X$ is given by the pre-composition $H_f(\alpha):= \alpha\circ f$. If we assume that $H$ has a bottom element $\bot$, the above doctrine is elementary taking as fibered equality the function $\delta_X: X\times X\to H$ such that $\delta_X(x_1,x_2)$ is equal to $\top$ if $x_1=x_2$, otherwise it is equal to $\bot$.
\item\label{example:ed-of-mltt} Consider the intensional version of the Martin-L\"{o}f intuitionistic type theory with the usual type constructors and a universe (of \textit{small types}), such as the one described in \cite{nordstrom1990programming}. For such theory, we obtain the syntactic category $\ml$ whose objects are small closed types, i.e.\ types in the empty context, and arrows are terms up to \textit{function extensionality}: two terms $x:X\vdash t(x):Y$ and $x:X\vdash t'(x):Y$ are equivalent if there is a term of the type
\[	x:X\vdash \mathsf{Id}_Y(t(x), t'(x)).\]
The composition of arrows is given by the substitution of terms.
This category has strict finite products and weak pullbacks due to the intensional identity types. Consider the functor
\begin{equation*}
	F^{ML}:\ml\op\to\mathsf{InfSL}
\end{equation*}
which sends a closed type $X$ to $F^{ML}(X)$ defined as the poset of equivalence classes of types depending on $X$ with respect to \textit{equiprovability}: two types $B(x)$ and $B'(x)$ with $x:X$ are equivalent if there is a term of 
\[x:A\vdash (B(x)\to B'(x))\times(B'(x)\to B(x)).\]
The action of $F^{ML}$ on the arrows of $\ml$ is given by substitution. This functor is an elementary doctrine with fibered equality given by the identity types. The doctrine is also existential, universal, implicational with full weak comprehensions.
\item Consider the Minimalist Foundation ideated by Maietti and Sambin in \cite{maiettisambin005} and then fully formalized by Maietti in \cite{maietti2009minimalist}. This theory is a dependent type theory with an \textit{intensional level}, called minimal type theory (\textbf{mTT}), and an \textit{extensional level} (\textbf{emTT}). The intensional level gives rise to a syntactic category $\mathcal{CM}$ whose objects are closed \textit{collections} and arrows are terms between collections up to function extensionality. This category has strict finite products and weak pullbacks. We consider the doctrine 
\begin{equation}
G^{\textbf{mTT}}:\mathcal{CM}\op\to \mathsf{InfSL}
\end{equation}
which sends a closed collection $X$ to the equivalence class of  \textit{propositions} depending on $X$, i.e.\ $x:X\vdash B(x)\ prop$ up to equiprovability. This doctrine, similarly to the above example, is elementary, existential, universal, implicational with full weak comprehensions and comprehensive diagonals.
		\end{itemize}
\end{examples}

Primary doctrines form a 2-category denoted by \textbf{PD}. Objects of \textbf{PD} are primary doctrines and 1-arrows from $\mathsf{P}$ to $\mathsf{P}'$ are pairs $(F,f)$ where $F:\mathscr{C}\to \CC'$ is a functor which preserves finite products and $f:\PP\Rightarrow \PP'$ is a natural transformation such that, for every object $X\in \CC$, the functor $f_X:\mathsf{P}(X)\to \mathsf{P}'(F(X))$ preserves finite meets.

\begin{equation}\label{diagram: 1 arrow PD}
\begin{tikzcd}
	{\mathscr{C}\op} \\
	&&& {\mathsf{InfSL}} \\
	{{\mathscr{C}'}\op}
	\arrow["{F\op}"', from=1-1, to=3-1]
	\arrow[""{name=0, anchor=center, inner sep=0}, "{\mathsf{P}'}"', from=3-1, to=2-4]
	\arrow[""{name=1, anchor=center, inner sep=0}, "{\mathsf{P}}", from=1-1, to=2-4]
	\arrow["f"', shorten <=5pt, shorten >=5pt, from=1, to=0]
\end{tikzcd}\end{equation}
A 2-arrow $\theta:(F,f)\to (G,g)$ is a natural transformation $\theta:F\Rightarrow G$ as in the following diagram 

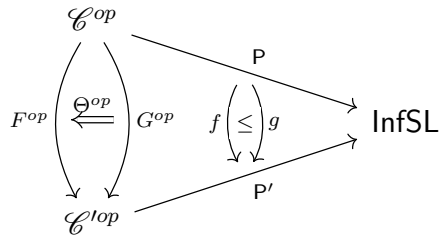
\begin{equation}\label{diagram: 2 arroew PD}\begin{tikzcd}
	{\mathscr{C}\op} \\
	&&& {\mathsf{InfSL}} \\
	{{\mathscr{C}'}\op}
	\arrow[""{name=0, anchor=center, inner sep=0}, "{F\op}"', bend right, from=1-1, to=3-1]
	\arrow[""{name=1, anchor=center, inner sep=0}, "{\mathsf{P}'}"', from=3-1, to=2-4]
	\arrow[""{name=2, anchor=center, inner sep=0}, "{\mathsf{P}}", from=1-1, to=2-4]
	\arrow[""{name=3, anchor=center, inner sep=0}, "{G\op}", bend left, from=1-1, to=3-1]
	\arrow[""{name=4, anchor=center, inner sep=0}, "g", bend left, shorten <=5pt, shorten >=5pt, from=2, to=1]
	\arrow[""{name=5, anchor=center, inner sep=0}, "f"', bend right, shorten <=5pt, shorten >=5pt, from=2, to=1]
	\arrow["{\Theta\op}"', shorten <=5pt, shorten >=5pt, Rightarrow, from=3, to=0]
	\arrow["\le"{description}, draw=none, from=4, to=5]
\end{tikzcd}\end{equation}
such that $f_X(\alpha)\le \PP'_{\Theta_X}(g_X(\alpha))$ for every $X\in\CC$ and $\alpha\in\PP(X)$.

We denote with \textbf{ED} the 2-full 2-subcategory of \textbf{PD} given by elementary doctrines and  1-arrows $(F,f)$ of \textbf{PD} such that the fibered equality is preserved, that is: for every $X\in\CC$
\begin{equation*}
	f_{X\times X}(\delta_X)=\PP'_{\langle F(pr_1),F(pr_2)\rangle}(\delta_{F(X)}).
\end{equation*}
In \cite{pasqualicofree} the author provides a 2-right adjoint to the forgetful 2-functor $\textbf{ED}\to \textbf{PD}$.

	We denote with \textbf{EqD} the 2-full 2-subcategory of \textbf{ED} given by elementary doctrines with full comprehensions and comprehensive diagonals and 1-arrows of \textbf{ED} preserving comprehensions.
	
	The structure of primary doctrines provides a general framework in which one can define a notion of equivalence relation and its quotient as follows. 
	
	\begin{dfn}\label{dfn:P-eq relation}\normalfont
		Let $\Pdoc$ be a primary doctrine and $X$ an object of $\CC$. A \textit{$\PP$-equivalence relation} on $X$ is an element $\rho\in\PP(X\times X)$ such that
		\begin{itemize}
			\item[]reflexivity: $\top_X\le \PP_{\Delta_X}\rho$
			\item[]symmetry: $\PP_{\ang{p_2,p_1}}\rho\le \rho$
			\item[]transitivity: $\PP_{\ang{p_1,p_2}}\rho\wedge \PP_{\ang{p_2,p_3}}\rho\le\PP_{\ang{p_1,p_3}}\rho$.
		\end{itemize}
	\end{dfn}
	
	Observe that if $\Pdoc$ is elementary, then the fibered equality $\delta_{X}$ is a $\PP$-equivalence relation for every object $X\in\CC$.
	
	\begin{dfn}\label{dfn:quotient}\normalfont
		\normalfont	Let $\Pdoc$ be an elementary doctrine and let $\rho$ be a $\PP$-equivalence relation on $X$ in $\CC$.
		A \textit{quotient} of $\rho$ is an arrow $q:X\to C$ such that $\rho\le \mathsf{P}_{q\times q}\delta_C$ and, for every arrow $g:X\to Z$ such that $\rho\le \mathsf{P}_{g\times g}\delta_Z$, there exists a unique arrow $h:C\to Z$ such that $g=h\circ q$. Such a quotient is said \textit{stable} when, for every arrow $f:C'\to C$ and every pullback square of $\CC$ of the form
	
	$$\begin{tikzcd}
		X'\arrow[]{r}{q'} \arrow[swap]{d}{f'} & C' \arrow[]{d}{f}\\
		X\arrow[swap]{r}{q}& C
	\end{tikzcd}$$
 we have that the arrow $q':X'\to C'$ is a quotient of a $\PP$-equivalence relation. If $f:A\to B$ is an arrow in $\CC$, the $\PP$-\textit{kernel} of $f$ is the $\PP$-equivalence relation $\mathsf{P}_{f\times f}(\delta_B)$. A quotient $q:X\to C$ of the $\PP$-equivalence relation $\rho$ is called \textit{effective} if its $\PP$-kernel is $\rho$. The quotient $q$ is of \textit{effective descent} if the functor $$\mathsf{P}_q:\mathsf{P}(C)\to \mathpzc{Des}({\rho})$$ is an isomorphism.
\end{dfn}

	We denote with \textbf{QD} the 2-full 2-subcategory of \textbf{EqD} whose objects are elementary doctrines $\Pdoc$ of \textbf{EqD} with stable effective quotients of $\PP$-equivalence relations of effective descent. The 1-arrows of \textbf{QD} are the 1-arrows of \textbf{EqD} which preserve quotients.
	
	We now recall the construction of \textit{elementary quotient completion}, which was introduced in \cite{maietti2013quotient}. Given an elementary doctrine $\Pdoc$ in \textbf{ED}, define the category $\overline{\CC}$ whose
	\begin{itemize}
		\item[-] {objects} are pairs $(X,\rho)$, where $X\in\CC$ and $\rho$ is a $\PP$-equivalence relation on $X$,
		\item[-] {arrows} between two objects $(X,\rho)$ and $(Y,\sigma)$ are equivalence classes of the arrows $f:X\to Y$ such that $\rho\le\nolinebreak \mathsf{P}_{f\times f}(\sigma)$. Two arrows $f,f'$ are equivalent when $\rho\le \mathsf{P}_{f\times f'}(\sigma)$.
	\end{itemize}
	Consider the functor
	\begin{equation}\label{equation: elementary quot comp}
		\overline{\mathsf{P}}:\overline{\CC}\op\to  \mathsf{InfSL}
	\end{equation}
	which sends an object $(X,\rho)\in\overline{\CC}$ to $\overline{\mathsf{P}}(X,\rho):=\mathpzc{Des}({\rho})$ and an  arrow $\lfloor f\rceil$ to $\overline{\mathsf{P}}_{\lfloor f\rceil}:=\mathsf{P}_f$.
	The functor $\overline{\PP}$ is called the \textit{elementary quotient completion} of $\PP$ and it has stable effective quotients of $\overline{\PP}$-equivalence relations, and quotients are of effective descent. $\overline{\PP}$ has also comprehensive diagonals, see \cite{maietti2013quotient,maietti2013elementary}.
	
	Moreover, there is a canonical 1-arrow $$(J,j):\PP\to \overline{\mathsf{P}}$$  given by the functor $J$, which sends an object $X\in\CC$ to $(X,\delta_X)\in\overline{\CC}$ and an arrow $f:X\to Y$ to the arrow 
	$\lrfloor{f}:(X,\delta_X)\to (Y,\delta_Y)$.
	For every object $X\in\CC$ the natural transformation $j$ is defined as $j_X:=1_{\mathsf{P}(X)}$.
	
	The {elementary quotient completion} has the following universal property which appears as \cite[Theorem 5.8]{maietti2013quotient}.

	\begin{thm}[Maietti and Rosolini]\label{thm:elementary-quotient-completion-universal-property} For every elementary doctrine $\Pdoc$ in \emph{\textbf{EqD}}, the assignment $\PP\to \overline{\mathsf{P}}$ gives a left biadjoint to the forgetful 2-functor \emph{$U:\textbf{QD}\to\textbf{EqD}$}, i.e. the pre-composition with the 1-arrow $(J,j)$		
		induces an equivalence of categories
	\emph{	$$-\circ(J,j):\textbf{QD}(\overline{\mathsf{P}},R)\cong \textbf{EqD}(\PP,UR)$$}
		for every $R$ in $\mathbf{QD}$.
		\qed
	\end{thm}
	
	\begin{rmk}\label{rmk:step quotient completion}\normalfont
		In \cite{maietti2013elementary}, the authors make a deep analysis of universal constructions for elementary doctrines identifying three primary ones. Given an elementary doctrine, one can separately freely add effective quotients, comprehensions or enforce comprehensive diagonals. Combining these three constructions the authors extend the result in \Cref{thm:elementary-quotient-completion-universal-property}, providing a left biadjoint to the forgetful 2-functor $U:\mathbf{QD}\to \mathbf{ED}$.
		
	\end{rmk}
 \begin{rmk}\normalfont
 The elementary quotient completion $\Pdocc{\overline{\PP}}{\overline{\CC}}$ of an elementary doctrine $\Pdoc$ in \textbf{EqD} does not necessarily  provide an exact category, see \ref{examples: eqc}. However, the category $\overline{\CC}$ is at least regular, as shown in \cite[Proposition 4.15]{maietti2013quotient}.
 	\end{rmk}

\begin{rmk}\normalfont

In \cite{carbonimagno,carboni1998regular}, the authors provide an explicit description of the exact completion of a category with weak finite limits $\CC$, using the notion of \textit{pseudo equivalence relation }, which consists of a pair of arrows of $\CC$
\begin{equation*}
	r_1,r_2:R\to X
\end{equation*}
satisfying reflexivity, symmetry and transitivity condition. The exact completion of $\CC$ is category $\CC_{ex}$ whose objects are pseudo equivalence relations, and arrows between two objects $r_1,r_2:R\to X$ and $s_1,s_2:S\to Y$ are given by pairs $(f,\tilde{f})$ of arrows of $\CC$, $f:X\to Y$ and $\tilde{f}:R\to S$, such that $fr_i=s_i\tilde{f}$, for $i=1,2$, up to the notion of \textit{half-homotopy}. In case $\CC$ has strict finite products, pseudo equivalence relations are in bijections with $\Psi_{\CC}$-equivalence relations. Indeed, a pseudo equivalence relation $r_1,r_2:R\to X$ induces an element 
\begin{equation}\label{equation:correspondence per P-eqrel}
	\lrfloor{r:R\to X\times X}\in\Psi_{\CC}(X\times X),
\end{equation}
where $r$ is the unique arrow induced by the universal property of products, such that $\lrfloor{r}$ is a $\Psi_{\CC}$-equivalence relation on $X$. Moreover, a quotient of $\lrfloor{r}$ is a coequalizer of the parallel arrows $r_1,r_2$. These facts imply the following result, see \cite{maietti2013quotient,maietti2013elementary}.
\end{rmk}
	
	\begin{thm}\label{thm:exactcoompletionaseqm}
		Let $\CC$ be a category with strict products and weak pullbacks. Then the elementary quotient completion of the elementary doctrine of weak subobjects $\Pdocc{\Psi_{\CC}}{\CC}$ is equivalent to the elementary doctrine $\Pdocc{\mathsf{Sub}_{\CC_{ex}}}{\CC_{ex}}$.
		\qed
	\end{thm}

	\begin{examples}\label{examples: eqc}\normalfont
			We end this section with two examples of elementary quotient completion arising from type theory. 
	\begin{itemize}
		\item The elementary quotient completion of the doctrine $F^{ML}:\ml\op\to\mathsf{InfSL}$ gives rise to a category $\overline{\textbf{ML}}$ which is equivalent to the category \std\ of \textit{total setoids} à la Bishop over Martin-L\"{o}f type theory \cite{bishop1967foundations,palmgren2005bishop}. However, the category \std\ is the exact completion of \ml\ and the doctrines $F^{ML}$ and $\Psi_{\ml}$ are equivalent, see \cite[\S 7.1]{maietti2013quotient}. See also \cite{van2018exact}, where the authors obtain the category of setoids through a different construction called \textit{homotopy exact completion}.
		\item The first example of elementary quotient completion that is not an exact completion is $\Pd{\overline{G^{\textbf{mTT}}}}{\overline{\mathcal{CM}}}$. The category $\overline{\mathcal{CM}}$ is not exact and doctrine $\overline{G^{\textbf{mTT}}}$ provides an algebraic account of the interpretation of \textbf{emTT} into \textbf{mTT}, which is referred to as the \textit{quotient model} in \cite{maietti2009minimalist}. Another example of elementary quotient completion which is not an instance of the exact completion is given by the category of \textit{homotopy setoids} studied in \cite[\S 2]{cioffophd}.
	\end{itemize}
	\end{examples}

\section{Biased elementary doctrines}\label{section:biased elementary doctrines}


In this section, we generalize the framework of elementary doctrines to encompass suitable functors from categories with weak finite products, rather than strict ones. As mentioned in the introduction, our main goal is to find a suitable extension of the elementary quotient completion to generalize the exact completion for categories with weak finite limits (and also weak finite products). To this end, we begin with a minimal weakening of the notion of doctrine, defined by contravariant functors from categories with weak finite products to the category of posets. We will refer to these functors as \textit{biased doctrines}.

Recall that a category $\CC$ has weak finite products if for every two objects $X,Y\in\CC$ there exists a diagram $\weakbinprodg{p}{X}{W}{Y}$ such that the following weakened universal property holds: for any two arrows $f:A\to X$ and $g:A\to Y$ there exists a arrow $h:A\to W$ (not necessarily unique) such that $\pp_1h=f$ and $\pp_2h=g$. The following are some examples of categories with weak finite products.

\begin{examples}\label{examples: w prod }\normalfont \phantom{a}
	\begin{enumerate}
	\item  The category $\set$ of sets has strict finite products. However, for every two sets $X$ and $Y$ the product of $X\times Y$ with any non-empty set $T$ is a weak product of $X$ and $Y$. For instance consider $$\binprod{X}{X\times Y \times \{0,1\}}{Y}.$$
	
\item\label{item: examples weak prod Ho(Top)} Let $\hotop$ be the category of topological spaces and continuous functions up to homotopy, which is equivalent to the category of fractions of $\mathbf{Top}$ modulo homotopy \cite{GabZis}. Homotopy pullbacks in $\mathbf{Top}$ provide weak pullbacks in $\hotop$. Hence, for every topological space $X$ the slice category $\hotop/X$ has weak finite products. A weak pullback of two (homotopy classes of) arrows $f:A\to X$ and $g:B\to X$ is given by the diagram $\weakbinprodg{\pi}{A}{W}{B}$,
where $$W:=\{(a,b,\gamma)\in A\times B\times X^I\vert\ \gamma(0)=f(a), \gamma(1)=g(b)\}$$ and $\pi_1, \pi_2$ are respectively the projections into the spaces $A$ and $B$. The presence of possibly non-homotopic paths in $X$ implies that $\pi_1,\pi_2$ are not necessarily jointly monic.

\item\label{item: examples weak prod G-sets} It is well-known that for a monad $(\mathbb{T},\epsilon,\mu)$ over the category of sets $\mathsf{Set}$ the Eilenberg-Moore category $\set^\mathbb{T}$ is exact and the Kleisli category $\set_\mathbb{T}$ is a projective cover of it. Hence, $\set^\mathbb{T}$ is the exact completion of $\set_\mathbb{T}$, see \cite{borceux2} and \cite{carboni1998regular}.\newline
In particular, if $G=(G,\cdot,e,(\text{-})^{-1})$ is a group then the Eilenberg-Moore category of the action monad is the topos of $G$-sets $\Gsets$. The Kleisli category of this monad is the category of free algebras $\setsG$ which is a projective cover of $\Gsets$ and we have $(\setsG)_{ex}\cong\Gsets$. The objects of $\setsG$ are sets and arrows $f:X\leadsto Y$ between sets are functions of the form $\langle f_1,f_2\rangle: X\to Y\times G$. The composition of two arrows $f:X\leadsto Y$ and $g:Y\leadsto Z$ is given by $\mu_Z(g\times 1_G)f$. The category $\setsG$ has weak terminal object and a weak products of two sets $X,Y$  given by 
\[
\begin{tikzcd}
	X & {X\times G\times Y \times G} & Y
	\arrow["{\langle p_3,p_4\rangle}", rightsquigarrow, from=1-2, to=1-3]
	\arrow["{\langle p_1,p_2\rangle}"', squiggly, from=1-2, to=1-1]
\end{tikzcd}\]
Since the composition of arrows of $\setsG$ involves the multiplication of elements of $G$, the projections $\langle p_1,p_2\rangle, \langle p_3,p_4\rangle$ are not necessarily jointly monic if the group $G$ is not trivial.

\item\label{item: examples weak prod Markov} Markov categories provide a categorical framework for synthetic probability theory. A Markov category $\CC$ is a symmetric monoidal category in which every object $X$  is equipped with a commutative comonoid structure given by a comultiplication $\mathsf{copy}_X:X\to X \otimes X$ and counit $\mathsf{del}_X: X\to I$ satisfying suitable commutativity, compatibility and naturality conditions. Examples of Markov categories are the category $\mathsf{FinStoch}$ of finite sets and stochastic matrices, and the category $\mathsf{Stoch}$ of measurable spaces and Markov kernels. In Markov categories the tensor product of two objects $X,Y$ is a weak product with the following projections $$X\cong X\otimes I\overset{1_X \otimes\mathsf{del}_Y}{\longleftarrow} X\otimes Y \overset{ \mathsf{del}_X\otimes 1_Y}{\longrightarrow} I\otimes Y \cong Y$$
Intuitively, the arrows of the form $I\to X$ correspond to distributions over $X$. Hence, since there could be different probability measures on $X \otimes Y$ with the same marginal distributions, the arrows $1_X \otimes\mathsf{del}_Y, \mathsf{del}_X\otimes 1_Y$ are not jointly monic. We refer to \cite{fritz2020synthetic} for an exhaustive treatment of Markov categories.
\end{enumerate}
\end{examples}


We now define the notion of \textit{biased primary doctrine} observing that the modular correspondence between the categorical structures and logic is especially beneficial; the conditions in  \Cref{dfn:primarydoctrine}, that provide the fibers to have finite meets and the reindexings to preserve them, are mutually independent and we can rewrite them when $\CC$ has just weak finite products as follows. 

\begin{dfn}\label{dfn:weak-primarydoctrine}\normalfont
	Let $\CC$ be a category with weak finite products. A \textit{biased primary doctrine}  is a functor $\PP:\CC\op\to \mathsf{Pos}$ which takes value in the category $\mathsf{InfSL}$ of inf-semilattices, i.e.: 
	\begin{enumerate}
		\item for every object $X\in\CC$, $\mathsf{P}(X)$ has finite meets
		\item for every arrow $f:X\to Y$ in $\CC$, the map $\mathsf{P}_f:\mathsf{P}(Y)\to \mathsf{P}(X)$ preserves finite meets.
	\end{enumerate}
\end{dfn}
	We denote with \textbf{BPD} the 2-category of biased primary doctrines with 1-arrows given by those pairs $(F,f)$ as in diagram (\ref{diagram: 1 arrow PD}) where the functor $F$ does not necessarily preserve weak finite products, and the natural transformation $f$ is such that $f_X$ preserves finite meets for every object $X\in\CC$. A 2-arrow $\theta:(F,f)\to (G,g)$ is given by natural transformation $\theta:F\Rightarrow G$ as in diagram (\ref{diagram: 2 arroew PD}) satisfying  $f_X(\alpha)\le \PP'_{\Theta_X}(g_X(\alpha))$ for every $X\in\CC$ and $\alpha\in\PP(X)$.

	One might find it surprising that in \textbf{BPD} we have considered functors that do not preserve weak finite products.  As we will see in the next sections, the adoption of this notion of 1-arrow is motivated by the fact that several morphisms of our interest between biased primary doctrines do not preserve weak finite products.

\begin{obs}\label{obs: BPD as fibrations}\normalfont
In \cite{maietti2013quotient}, the authors observe that the 2-category \textbf{PD} of primary doctrines is equivalent to a suitable 2-category of  Grothendieck fibrations. Namely, the faithful fibrations $p:\mathcal{E}\to \CC$ between categories with finite products such that $p$ preserves them (see also \cite{jacobs1999categorical}). Similarly, we obtain a correspondence between biased primary doctrines and faithful fibrations $p:\mathcal{E}\to \CC$ between categories with weak finite products such that $p$ preserves them. This correspondence extends to an equivalence between \textbf{BPD} and a suitable 2-category of the above faithful fibrations.
\end{obs}

Contrary to the definition of biased primary doctrine, the conditions in Proposition \ref{prop:equiv-elementarydoctrine}, which characterize the elementary structure, are interdependent and the (strict) products in the base category $\CC$ play a key role in it. Since two weak products of the same objects are not necessarily isomorphic, we have to devise a way so that fibered equalities, which shall now become \textit{biased fibered equalities}, interact appropriately. As we expect, an elementary doctrine shall satisfy also the following Definition \ref{dfn:weakelementarydoctrine} of \textit{biased elementary doctrine}, but we will discuss in detail the relationship between the two notions in this section.

Before proceeding, if $X$ is an object of a category $\CC$, the diagram $\weakbinprodg{p}{X}{W}{X}$ is a weak product and $\beta$ is an element of $\PP(W)$, we will adopt the same notation $\mathpzc{Des}({\beta})$ of \Cref{dfn:descentobjects} for the elements $\alpha\in \mathsf{P}(X)$ such that
\[ 	\mathsf{P}_{\mathsf{p}_1}(\alpha)\wedge\beta\le \mathsf{P}_{\mathsf{p}_2}(\alpha)\]
and we will refer to as the descent data of $\beta$.

	\begin{dfn}\label{dfn:weakelementarydoctrine}\normalfont
	Let $\CC$ be a category with weak finite products. A \textit{biased elementary doctrine} is a biased primary doctrine $\Pdoc$	such that, for every object $X\in\CC$ and for every weak product $\weakbinprodg{p}{X}{W}{X}$
	there exists an element $\bfeq{X}{p}\in \mathsf{P}(W)$ satisfying:
	\begin{enumerate}
		\item\label{item:weakelementarydoctrine-1} For every commutative diagram
		$\begin{tikzcd}[row sep= tiny]
			&& X \\
			X & W \\
			&& X
			\arrow["d", from=2-1, to=2-2]
			\arrow["{\mathsf{p}_1}"', from=2-2, to=1-3]
			\arrow["{\mathsf{p_2}}", from=2-2, to=3-3]
			\arrow["1_X"', bend right=20, from=2-1, to=3-3]
			\arrow["1_X", bend left=20, from=2-1, to=1-3]
		\end{tikzcd}$
		we have $\top_X\le \mathsf{P}_d\bfeq{X}{p}$.
		\item\label{item:weakelementarydoctrine-2}
		$\mathsf{P}(X)= {\mathpzc{Des}}({\bfeq{X}{p}})$, i.e. for every $\alpha\in\PP(X)$ 
		\[\PP_{\mathsf{p}_1}\alpha \wedge \bfeq{X}{p}\le\PP_{\mathsf{p}_2}\alpha.\]
		\item\label{item:weakelementarydoctrine-3'} For any weak product $\weakbinprodgp{p}{X'}{W'}{X'}$ and for every commutative diagram 
		\[\begin{tikzcd}[row sep=tiny, column sep=small]
			& {X'} && X \\
			{W'} && W \\
			& {X'} && X
			\arrow["{\mathsf{p}'_1}", from=2-1, to=1-2]
			\arrow["{\mathsf{p}'_2}"', from=2-1, to=3-2]
			\arrow["{\mathsf{p}_1}", from=2-3, to=1-4]
			\arrow["{\mathsf{p}_2}"', from=2-3, to=3-4]
			\arrow["f", from=1-2, to=1-4]
			\arrow["f"', from=3-2, to=3-4]
			\arrow["g", from=2-1, to=2-3]
		\end{tikzcd}\]
		we have $\bfeqp{X'}{p}\le \mathsf{P}_{g}\bfeq{X}{p}.	$
		\item\label{item:weakelementarydoctrine-4} For every commutative diagram
		
		\begin{equation}\label{diagram cond 4 biasedl elementary}
			\begin{tikzcd}[row sep= tiny]
			&&&&& X \\
			&& W \\
			&& W &&& X \\
			U \\
			&& W &&& X \\
			&& W \\
			&&&&& X
			\arrow["{\mathsf{r_1}}", from=4-1, to=2-3]
			\arrow["{t}"'{pos=0.65},from=4-1, to=3-3]
			\arrow["{t'}"{pos=0.65},from=4-1, to=5-3]
			\arrow["{\mathsf{r}_2}"', from=4-1, to=6-3]
			\arrow["{\mathsf{p}_2}"{pos=0.75}, from=5-3, to=7-6]
			\arrow["{\mathsf{p}_1}"{pos=0.6}, from=6-3, to=5-6, crossing over]
			\arrow["{\mathsf{p}_2}"'{pos=0.6}, from=6-3, to=7-6]
			\arrow["{\mathsf{p}_1}"{pos=0.6}, from=2-3, to=1-6]
			\arrow["{\mathsf{p}_1}"'{pos=0.75}, from=3-3, to=1-6]
			\arrow["{\mathsf{p}_2}"'{pos=0.6}, from=2-3, to=3-6, crossing over]
			\arrow["{\mathsf{p}_1}"{pos=0.75}, from=5-3, to=3-6]
			\arrow["{\mathsf{p}_2}"'{pos=0.75}, from=3-3, to=5-6, crossing over]
		\end{tikzcd}
	\end{equation}
		where $U$ is a weak product $\weakbinprodg{r}{W}{U}{W}$, we have $\bfeq{X}{p}\in {\mathpzc{Des}}({\mathsf{P}_{t}\bfeq{X}{p}\wedge \mathsf{P}_{t'}\bfeq{X}{p}})$, i.e.
		\[\PP_{\mathsf{r}_1}\bfeq{X}{p}\wedge\mathsf{P}_{t}\bfeq{X}{p}\wedge \mathsf{P}_{t'}\bfeq{X}{p}\le \PP_{\mathsf{r}_2}\bfeq{X}{p}. \]
	\end{enumerate}
\end{dfn}

	Now we discuss the apparent similarities between the above conditions and those in Proposition \ref{prop:equiv-elementarydoctrine}. Conditions \ref{item:weakelementarydoctrine-1} and \ref{item:weakelementarydoctrine-2} of Definition \ref{dfn:weakelementarydoctrine} are transcriptions of condition \ref{item:equiv-elementarydoctrine-1} and \ref{item:equiv-elementarydoctrine-2} of Proposition \ref{prop:equiv-elementarydoctrine} using weak binary products instead of strict ones. Condition \ref{item:weakelementarydoctrine-3'} of \Cref{dfn:weakelementarydoctrine} holds in the context of elementary doctrines. It follows directly from condition \ref{item:elementarydoctrine-1} of Definition \ref{dfn:elementarydoctrine}, which can be in turn derived from the conditions \ref{item:equiv-elementarydoctrine-1}, \ref{item:equiv-elementarydoctrine-2} and \ref{item:equiv-elementarydoctrine-3} of Proposition \ref{prop:equiv-elementarydoctrine}.
Also condition \ref{item:weakelementarydoctrine-4} of Definition \ref{dfn:weakelementarydoctrine} holds in the context of elementary doctrines. Indeed, in case of strict finite products, the inequality in condition \ref{item:weakelementarydoctrine-4} of Definition \ref{dfn:weakelementarydoctrine} becomes
\[\PP_{\mathsf{r}_1}\delta_X\wedge{\mathsf{P}_{\ang{p_1,p_3}}\delta_X\wedge \mathsf{P}_{\ang{p_2,p_4}}\delta_X}\le \PP_{\mathsf{r}_2}\delta_X\]
which holds thanks to conditions \ref{item:equiv-elementarydoctrine-2} and \ref{item:equiv-elementarydoctrine-3} of \Cref{prop:equiv-elementarydoctrine}. The main reason why we explicitly assume conditions \ref{item:weakelementarydoctrine-3'} and \ref{item:weakelementarydoctrine-4} in \Cref{dfn:weakelementarydoctrine} is that our main examples of biased elementary doctrines lack of a condition analogous to \ref{item:equiv-elementarydoctrine-3} of \Cref{prop:equiv-elementarydoctrine}, which is necessary to derive them.

 For similar reasons, we are not able to provide an equivalent formulation of biased elementary doctrines in style of \Cref{dfn:elementarydoctrine} using left adjoints along parameterized diagonals, since it would imply a condition similar to \ref{item:equiv-elementarydoctrine-3} of \Cref{prop:equiv-elementarydoctrine}.

 Moreover, a peculiarity of (strict) elementary doctrines is the equality of the two elements 
\begin{equation}\label{equation:delta of prod = prod of deltas}
	\delta_{X\times Y}=\delta_{X} \boxtimes \delta_{Y}.
\end{equation}
which means that the equality on $X\times Y$ coincides with the component-wise equality on $X$ and $Y$.

 In the following \Cref{lemma-conditions-3-4} we prove that the inequality $\delta_{X\times Y}\le\delta_{X} \boxtimes  \delta_{Y}$, restated using weak finite products, still holds for the biased elementary doctrines. As already mentioned, the opposite inequality $\delta_{X} \boxtimes  \delta_{Y}\le \delta_{X\times Y}$ (which is condition \ref{item:equiv-elementarydoctrine-3} of \Cref{prop:equiv-elementarydoctrine}) does not necessarily hold in the context of biased elementary doctrines, as we will see in \Cref{example: toy difference of equalities,example: why-proof-irrelevant,example:slicedoctrine-of-ed,example: proof-irrelevant g-sets}.

\begin{lemma}\label{lemma-conditions-3-4}
	If $\Pdoc$ is a biased elementary doctrine, then the following conditions hold
	\begin{itemize}
		\item\label{item:weakelementarydoctrine-3} Given two objects $X,X'\in\CC$, then for any weak products $\weakbinprodg{q}{X}{V}{X'}$, $\weakbinprodgp{p}{X}{W}{X}$, $\weakbinprodgp{p}{X'}{W'}{X'}$, and 
		$\weakbinprodg{r}{V}{U}{V}$ and for every commutative diagram
		
		\begin{equation}\label{diagram: first diagram proof-relevant / proof - irrelevant equlity}
		\begin{tikzcd}[row sep= tiny]
			&&&&& X \\
			&& V \\
			&& W &&& X' \\
			U \\
			&& W' &&& X \\
			&& V \\
			&&&&& X'
			\arrow["{\mathsf{r_1}}", from=4-1, to=2-3]
			\arrow["{t}"'{pos=0.65},from=4-1, to=3-3]
			\arrow["{t'}"{pos=0.65},from=4-1, to=5-3]
			\arrow["{\mathsf{r}_2}"', from=4-1, to=6-3]
			\arrow["{\mathsf{p}'_2}"{pos=0.75}, from=5-3, to=7-6]
			\arrow["{\mathsf{q}_1}"{pos=0.6}, from=6-3, to=5-6, crossing over]
			\arrow["{\mathsf{q}_2}"'{pos=0.6}, from=6-3, to=7-6]
			\arrow["{\mathsf{q}_1}"{pos=0.6}, from=2-3, to=1-6]
			\arrow["{\mathsf{p}_1}"'{pos=0.75}, from=3-3, to=1-6]
			\arrow["{\mathsf{q}_2}"'{pos=0.6}, from=2-3, to=3-6, crossing over]
			\arrow["{\mathsf{p}'_1}"{pos=0.75}, from=5-3, to=3-6]
			\arrow["{\mathsf{p}_2}"'{pos=0.75}, from=3-3, to=5-6, crossing over]
		\end{tikzcd}
	\end{equation}
		it follows that
		\[	\bfeq{V}{r}\le \mathsf{P}_{t}\bfeq{X}{p} \wedge \mathsf{P}_{t'}\bfeqp{X'}{p}.\] 
		\item\label{item:weakelementarydoctrine-5} If $\weakbinprodg{p}{X}{W}{X}$ and $\weakbinprodgp{p}{X}{W'}{X}$ are two weak products of $X$ and $X$, then for any commutative diagram 
		$$\begin{tikzcd}[row sep= tiny]
			&& X \\
			W' & W \\
			&& X
			\arrow["h", from=2-1, to=2-2]
			\arrow["{\mathsf{p}_1}"', from=2-2, to=1-3]
			\arrow["{\mathsf{p_2}}", from=2-2, to=3-3]
			\arrow["\mathsf{p}'_2"', bend right=20, from=2-1, to=3-3]
			\arrow["\mathsf{p}'_1", bend left=20, from=2-1, to=1-3]
		\end{tikzcd}$$
		it follows that  $\bfeqp{X}{p}\le \mathsf{P}_h(\bfeq{X}{p}).$
	\end{itemize}
	\begin{proof}\normalfont
		The two results follow directly from condition \ref{item:weakelementarydoctrine-3'} of Definition \ref{dfn:weakelementarydoctrine}.
\end{proof}
\end{lemma}

\begin{rmk}\label{rmk: proof-relevant/irrelevant equality}\normalfont
In the situation of the above lemma, considering a diagram as (\ref{diagram: first diagram proof-relevant / proof - irrelevant equlity}) we will refer to  $\bfeq{V}{r}$ as the \textit{proof-relevant equality} of $V$ and to $\mathsf{P}_{t}\bfeq{X}{p} \wedge \mathsf{P}_{t'}\bfeqp{X'}{p}$ as the \textit{proof-irrelevant} or \textit{component-wise equality} of $V$. \Cref{lemma-conditions-3-4} states that proof-relevant equality implies proof-irrelevant one, the vice versa does not necessarily hold as shown in the following example.
\end{rmk}
\begin{example}\label{example: toy difference of equalities}\normalfont
A trivial example of the fact that the proof-relevant and the proof-irrelevant equalities may non coincide is given by the doctrine of subobjects $\mathsf{Sub}_\set$ of $\set$. Indeed, given two sets $X,Y$ and a weak product  $V:={X\times Y \times \{0,1\}}$ we can consider the (weak) product $U:= V\times V$. Proof-relevant equality is given by the equality in $V$ while the proof-irrelevant one is given by the equalities of $X$ and $Y$ components. For instance, the elements $(x,y,0)$ and $(x,y,1)$ of $V$ are proof-irrelevant equal but not proof-relevant equal.
\end{example}
We now define a crucial property for elements in the fibers of weak finite products.
\begin{dfn}\label{dfn:reindex-determined-by-projections}\normalfont
	Let $\Pdoc$ be a biased primary doctrine and let $X_1,\dots,X_n$ be objects of $\CC$. Given a weak product $\pp_i:W\to X_i$, for $\lel{1}{i}{n}$, of $X_1,\dots,X_n$, an element $\beta\in \mathsf{P}(W)$ is called \textit{reindexed by projections (rbp)} if for any two arrows $h,h':Z\to W$, such that $\mathsf{p}_i\circ h=\mathsf{p}_i\circ h'$, for $\lel{1}{i}{n}$, we have \[\mathsf{P}_{h}(\beta)=\mathsf{P}_{h'}(\beta).\]
\end{dfn}
 A first example of element that is rbp is given by the biased fibered equality of Definition \ref{dfn:weakelementarydoctrine}. 

	\begin{lemma}\label{lemma:delta-is-determined-by-projections}
	If $\Pdoc$ is a biased elementary doctrine and  $\weakbinprodg{p}{X}{W}{X}$ is a weak product in $\CC$, then $\bfeq{X}{p}$ is reindexed by projections.
	\begin{proof}\normalfont
		 Let $h_1,h_2:Z\to W$ be two arrows such that $\mathsf{p}_i\circ h_1=\mathsf{p}_i\circ h_2$, for $i=1,2$, and consider a commutative diagram of the form 
		 \[\begin{tikzcd}[row sep= small, column sep =small]
		 	&& Z && W \\
		 	Z & V && U \\
		 	&& Z && W
		 	\arrow["{\mathsf{q}_1}", from=2-2, to=1-3]
		 	\arrow["{\mathsf{q}_2}"', from=2-2, to=3-3]
		 	\arrow["{\mathsf{r}_1}", from=2-4, to=1-5]
		 	\arrow["{\mathsf{r}_2}"', from=2-4, to=3-5]
		 	\arrow["{h_1}", from=1-3, to=1-5]
		 	\arrow["{h_2}"', from=3-3, to=3-5]
		 	\arrow["g", from=2-2, to=2-4]
		 	\arrow["d", from=2-1, to=2-2]
		 	\arrow["1_Z", bend left= 15, from=2-1, to=1-3]
		 	\arrow[swap,"1_Z",bend right=15, from=2-1, to=3-3]
		 \end{tikzcd}\]
		where $\weakbinprodg{q}{Z}{V}{Z}$ and $\weakbinprodg{r}{W}{U}{W}$ are weak products. Conditions \ref{item:weakelementarydoctrine-1} and \ref{item:weakelementarydoctrine-3'} of Definition \ref{dfn:weakelementarydoctrine} imply that $\PP_d\mathsf{P}_{g}\mathsf{P}_{t}\bfeq{X}{p}=\top_Z=\PP_d\mathsf{P}_{g}\mathsf{P}_{t'}\bfeq{X}{p}$, where $t,t':U\to W$ are two arrows such as in the diagram (\ref{diagram cond 4 biasedl elementary}). The inequality $\mathsf{P}_{h_1}\bfeq{X}{p}\le \mathsf{P}_{h_2}\bfeq{X}{p}$ is obtained as follows:
		\begin{align}
			\mathsf{P}_{h_1}\bfeq{X}{p}&=\PP_d\mathsf{P}_{g}\mathsf{P}_{{\mathsf{r}_1}}\bfeq{X}{p} \wedge \top_Z \tag*{} \\
			&=\PP_d\mathsf{P}_{g}(\mathsf{P}_{{\mathsf{r}_1}}\bfeq{X}{p}\wedge \mathsf{P}_{t}\bfeq{X}{p}\wedge \mathsf{P}_{t'}\bfeq{X}{p}) \tag*{}\\ 
			&\le \PP_d\mathsf{P}_{g}\mathsf{P}_{{\mathsf{r}_2}}\bfeq{X}{p}\tag*{}\\ 
			&=\mathsf{P}_{h_2}\bfeq{X}{p}. \tag*{}
		\end{align}
		The opposite inequality is obtained similarly considering an arrow $g':Z\to U$ such that ${\mathsf{r}_2}g'=h_1$ and ${\mathsf{r}_1}g'=h_2$.
	\end{proof}
\end{lemma}
The above lemma implies a refinement of the second condition of \Cref{lemma-conditions-3-4}.
\begin{cor}\label{cor:delta-reindexings}
	Let $\Pdoc$ be a biased elementary doctrine and let $X$ be an object of $\CC$. If $\weakbinprodg{p}{X}{W}{X}$ and $\weakbinprodgp{p}{X}{W'}{X}$ are two weak products and $h:W'\to W$ is an arrow satisfying $\mathsf{p}_ih=\mathsf{p}'_i$ for $i=1,2$, then $$\mathsf{P}_h\bfeq{X}{p}= \bfeqp{X}{p}.$$
	\qed
\end{cor}

We denote with \textbf{BED} the 2-full 2-subcategory of \textbf{BPD} given by biased elementary doctrines and 1-arrows $(F,f)$ of \textbf{BPD} such that, for every object $X\in\CC$ and for every weak product $\weakbinprodg{p}{X}{W}{X}$  in the domain of $F$ we have
\[f_W(\bfeq{W}{p})= \PP'_h(\bfeq{V}{q})\]

for any commutative diagram of the form
\begin{equation}\label{diagram: BED}
\begin{tikzcd}
	& V \\
	{F(X)} && {F(X)} \\
	& {F(W)}
	\arrow["{\mathsf{q}_1}"', from=1-2, to=2-1]
	\arrow["{\mathsf{q}_2}", from=1-2, to=2-3]
	\arrow["{F(\mathsf{p}_1)}", from=3-2, to=2-1]
	\arrow["{F(\mathsf{p}_2)}"', from=3-2, to=2-3]
	\arrow["h", dashed, from=3-2, to=1-2]
\end{tikzcd}\end{equation}
where  $\weakbinprodg{q}{F(X)}{V}{F(X)}$ is a weak product. Since we do not assume 1-arrows  between biased elementary doctrines to preserve weak finite products, the above condition represents a mild requirement for preservation of biased fibered equalities. Moreover, thanks to \Cref{lemma:delta-is-determined-by-projections} and \Cref{cor:delta-reindexings} it is enough that, for every weak product $\weakbinprodg{p}{X}{W}{X}$, the above condition is satisfied for just one diagram as in (\ref{diagram: BED}).

Definition \ref{dfn:weakelementarydoctrine} requires, for any object $X\in \CC$, the verification of certain conditions for every choice of weak product $\weakbinprodg{p}{X}{W}{X}$. We say that a category $\CC$ has\textit{ a choice of weak products} if there exists a functor $\omega:\CC\times\CC\to \CC^\Lambda$, where $\Lambda$ is the category with three objects and the non trivial span 
\[\begin{tikzcd}
	\bullet & \bullet & \bullet
	\arrow[from=1-2, to=1-1]
	\arrow[from=1-2, to=1-3]
\end{tikzcd},\] such that:
\begin{itemize}
	\item[-] the value $\omega(X,Y)$ is a weak product of the objects $X,Y\in\CC$, that we shall denote with the usual notation $\binprod{X}{X\times Y}{Y}$,
	\item[-] the arrow $\omega(f,g):\omega(X,Y)\to \omega(A,B)$ satisfies $p_1\circ\omega(f,g)=f\circ {p}_1$ and $p_2\circ\omega(f,g)=g\circ {p}_2$, for every pair of arrows $f:X\to A$ and $g:Y\to B$ of $\CC$. We will adopt the usual notation $f\times g$ for the choice of the weak product of arrows $\omega(f,g)$.
\end{itemize}
We now prove that a choice of weak products provides an easier description of the biased elementary doctrines.

\begin{thm}\label{thm:equiv-weakelementarydoctrines}
	If $\Pdoc$ is a biased primary doctrine such that $\CC$ has a choice of weak products and such that for every object $X\in\CC$ there exists an element $\delta_X\in \mathsf{P}(X\times X)$ satisfying:
	\begin{enumerate}
		\item\label{item:equiv-weakelementarydoctrines-1} there exists an arrow $X\overset{d}{\longrightarrow}X\times X$ such that $1_X=p_1\circ d=p_2\circ d$ and $\top_X\le \mathsf{P}_d(\delta_X)$,
		\item\label{item:equiv-weakelementarydoctrines-2}
		$\mathsf{P}(X)= {\mathpzc{Des}}({\delta_X})$,
		\item\label{item:equiv-weakelementarydoctrines-3} If $f:Y\to X$ is an arrow of $\CC$, then	$\delta_{Y}\le \mathsf{P}_{f\times f}\delta_{X}$,
		\item\label{item:equiv-weakelementarydoctrines-4} 
			$\delta_X\in {\mathpzc{Des}}({\mathsf{P}_{{p_1 \times p_1}}\delta_X\wedge \mathsf{P}_{{p_2 \times p_2} }\delta_X})$,
		\end{enumerate}
		then the functor $\Pdoc$ is a biased elementary doctrine.
		\begin{proof}\normalfont
			First observe that the hypothesis of the theorem imply that $\delta_X$ is rbp. The proof is the same of \Cref{lemma:delta-is-determined-by-projections}.
			Hence, for each weak product $\weakbinprodg{p}{X}{W}{X}$, the weak universal property of weak products induces an arrow $h:W\to X\times X$ such that $p_i\circ h=\mathsf{p}_i$, for $i=1,2$. Even if $h$ is not unique, we can uniquely reindex $\delta_X$ along such arrows and define $\bfeq{X}{p}:=\mathsf{P}_h\delta_{X}$, which trivially satisfies conditions \ref{item:weakelementarydoctrine-1} and \ref{item:weakelementarydoctrine-2} of Definition \ref{dfn:weakelementarydoctrine}. We now prove condition \ref{item:weakelementarydoctrine-3'} of Definition \ref{dfn:weakelementarydoctrine}. Let $\weakbinprodg{p}{X}{W}{X}$, $\weakbinprodgp{p}{X'}{W}{X'}$ be two weak products and let  $f:X'\to X$ and $g:W'\to W$ be two arrows making the diagram in condition \ref{item:weakelementarydoctrine-3'} commute. The weak universal property of weak products induces two arrows $h:W\to X\times X$ and $k:W'\to X'\times X'$ such that $p_i\circ h=\mathsf{p}_i$ and $p_i\circ k=\mathsf{p}_i$ , for $i=1,2$. We obtain the inequality $ \bfeqp{X'}{p}\le \mathsf{P}_{g}\bfeq{X}{p}$ as follows
			\begin{align}
				\bfeqp{X'}{p}&:= \mathsf{P}_k\delta_{X'} \tag*{} \\
				&\le \mathsf{P}_k\mathsf{P}_{f\times f}\delta_{X} \tag*{}\\
				&= \mathsf{P}_g\mathsf{P}_h\delta_{X} \tag*{}\\
				&:=\mathsf{P}_{g}\bfeq{X}{p}. \tag*{}
			\end{align} \\
			We now prove condition \ref{item:weakelementarydoctrine-4} of Definition \ref{dfn:weakelementarydoctrine}. Let $\weakbinprodg{r}{W}{U}{W}$ be a weak product and let  $t,t':U\to W$ be two arrows as in the commutative diagram of condition \ref{item:weakelementarydoctrine-4}.
			The weak universal property of weak products induces arrows $h:W\to X\times X$ and $k:U \to (X\times X)\times(X\times X)$ satisfying $p_i\circ h=\mathsf{p}_i$ and  $p_i\circ k=h\circ \mathsf{p}_i$, for $i=1,2$. We obtain the relation $\bfeq{X}{p}\in {\mathpzc{Des}}({\mathsf{P}_{t}\bfeq{X}{p}\wedge \mathsf{P}_{t'}\bfeq{X}{p}})$ as follows:
			\begin{align}
				&\mathsf{P}_{\mathsf{p}_1}\bfeq{X}{p} \wedge \mathsf{P}_{t}\bfeq{X}{p}\wedge \mathsf{P}_{t'}\bfeq{X}{p}\tag*{}\\
				&:= \mathsf{P}_{\mathsf{p}_1}\mathsf{P}_h\delta_{X}\wedge \mathsf{P}_{t}\mathsf{P}_h\delta_X\wedge \mathsf{P}_{t'}\mathsf{P}_h\delta_X	\tag*{} \\
				&=\mathsf{P}_k\mathsf{P}_{p_1}\delta_{X}\wedge \mathsf{P}_k\mathsf{P}_{p_1 \times p_1}\delta_{X} \wedge \mathsf{P}_k\mathsf{P}_{p_2\times p_2}\delta_{X} \tag*{}\\
				&= \mathsf{P}_k(\mathsf{P}_{p_1}\delta_{X}\wedge \mathsf{P}_{p_1\times p_1}\delta_{X} \wedge \mathsf{P}_{p_2\times p_2}\delta_{X}) \tag*{}\\
				&\le \mathsf{P}_k\mathsf{P}_{p_2}\delta_{X} \tag*{} \\
				&=\mathsf{P}_{\mathsf{p}_2}\mathsf{P}_h\delta_{X} \tag*{}\\
				&:=\mathsf{P}_{\mathsf{p}_2}\bfeq{X}{p}. \tag*{}
			\end{align}
		\end{proof}
	\end{thm}

In the contexts of biased elementary doctrines, (full) (weak) comprehensions are defined as in Definition \ref{dfn:ed-comprehensions}. Instead, comprehensive diagonals are defined as follows.

\begin{dfn}\label{dfn:comprhe. diagonals BED}\normalfont
	A biased elementary doctrine $\Pdoc$ has \textit{comprehensive diagonals} if for every pair of arrows $f,g:A\to X$ then $f=g$ if and only if there exists an arrow $h$ in a commutative diagram of the form
	\[\begin{tikzcd}[row sep= tiny]
		&& X \\
		A & W \\
		&& X
		\arrow["h", from=2-1, to=2-2]
		\arrow["{\mathsf{p}_1}"', from=2-2, to=1-3]
		\arrow["{\mathsf{p_2}}", from=2-2, to=3-3]
		\arrow["g"', bend right=20, from=2-1, to=3-3]
		\arrow["f", bend left=20, from=2-1, to=1-3]
	\end{tikzcd}\]
	where $\weakbinprodg{p}{X}{W}{X}$ is a weak product, such that $\top_A\le\PP_h\bfeq{X}{p}$.
\end{dfn}
Thanks to \Cref{lemma:delta-is-determined-by-projections}, in the above definition the existence of such an arrow $h$ with $\top_A\le\PP_h\bfeq{X}{p}$ implies that for every $h'$ making the diagram commute $\top_A\le\PP_{h'}\bfeq{X}{p}$.

\example\label{example:wed-of-sed}\normalfont
Every elementary doctrine is actually a biased elementary doctrine. Indeed, if $\Pdoc$ is elementary, then for every $X\in\CC$ and weak product $\wbinprod{X}{W}{X}$, there exists a unique arrow $\langle \mathsf{p}_1,\mathsf{p}_2\rangle:W\to X\times X$ into the strict product $X\times X$ making the obvious diagram commute. The biased fibered equalities are given by the reindexings $\bfeq{X}{\pp}:=\mathsf{P}_{\langle \mathsf{p}_1,\mathsf{p}_2\rangle }(\delta_{X})\in \mathsf{P}(W)$. A simple verification shows that the conditions in Proposition \ref{prop:equiv-elementarydoctrine}  for $\delta_X$ imply that element $\bfeq{X}{\pp}$ satisfies the conditions in Definition \ref{dfn:weakelementarydoctrine}.
\endexample

\example\label{example: bed as precomposition}
\normalfont Let $\Pdoc$ be an elementary doctrine and let $F:\DD\to \CC$ be a functor from a category $\DD$ with weak finite products. A trivial computation shows that the composition $\PP\circ F\op$ is a biased elementary doctrine taking as biased fibered equality $\bfeq{X}{\pp}:= \PP_{\langle F(\pp_1), F(\pp_2)\rangle}\delta_{F(X)}$, for every $X\in\DD$ and weak product $\wbinprod{X}{W}{X}$. For instance, we can consider the category $\set_G$ of \Cref{examples: w prod } and the functor $R:\set_G\to \set$ which is the right adjoint of the adjunction generating the action monad. The composition of $R$ with some of the elementary doctrines in \Cref{examples: elementary doctrines}, provides examples of biased elementary doctrines such as $\mathsf{Sub}_\set R\op, \Psi_\set R\op$ and $HR\op$. The same holds if we assume $\PP$ biased. The following is another instance of this construction.
\endexample

\example{(Slice doctrines)}\label{example:slicedoctrine-of-ed} \normalfont Let $\Pdoc$ be an elementary doctrine such that $\CC$ has weak pullbacks. Hence, the slices $\CC/A$ have weak finite products for every object $A\in\CC$. We can consider an obvious functor $$\Pd{\mathsf{P}_{/A}}{(\CC/A)}$$ which sends an object $(f:X\to A)\in\CC/A$ into $\mathsf{P}_{/A}(f):=\mathsf{P}(X)$ and  an arrow $h:f\to g$ of $\CC/A$ as $\mathsf{P}_{/A}(h):=\mathsf{P}_h$. We will refer to the functor $\mathsf{P}_{/A}$ as the \textit{slice doctrine} over $A$. We claim that the slice doctrine $\Pd{\mathsf{P}_{/A}}{(\CC/A)}$ is a biased elementary doctrine for all object $A\in\CC$. Indeed, let $x:X\to A$ be an object of $\CC/A$ and consider a weak product of $\weakbinprodg{\pi}{x}{w}{x}$ given by the common value of the two composites of a weak pullback diagram
\begin{center}
	\begin{tikzcd}
		V\arrow[swap]{d}{\pi_1} \arrow[]{r}{\pi_2}  & X \arrow[]{d}{x} \\
		X \arrow[swap]{r}{x} & A.
	\end{tikzcd}
\end{center} 
The biased fibered equality is $\bfeq{x}{\pi}:=\mathsf{P}_{\ang{\pi_1,\pi_2}}\delta_{X}\in \mathsf{P}_{/A}(w) (\in\mathsf{P}(V))$. Conditions \ref{item:weakelementarydoctrine-1} and \ref{item:weakelementarydoctrine-2} of Definition \ref{dfn:weakelementarydoctrine} trivially follow from conditions  \ref{item:equiv-elementarydoctrine-1} and \ref{item:equiv-elementarydoctrine-2} of Proposition \ref{prop:equiv-elementarydoctrine}. We now prove that condition \ref{item:weakelementarydoctrine-3'} of Definition \ref{dfn:weakelementarydoctrine} holds. Let $x':X'\to A$ be an object of $\CC/A$ and let $w':V' \to A$ be a weak product $\weakbinprodgp{\pi}{x'}{w'}{x'}$. Given two arrows $f:x'\to x$ and $g:w'\to w$ such that $\pi_i\circ g= f\circ \pi_i$, for $i=1,2$, then we obtain condition \ref{item:weakelementarydoctrine-3'} as follows
\begin{align}
	(\mathsf{P}_{/A})_g\bfeq{x}{\pi}&:=\mathsf{P}_g\mathsf{P}_{\ang{\pi_1,\pi_2}}\delta_{X}  \tag*{}\\
	&=\mathsf{P}_{\ang{\pi'_1,\pi'_2}}\mathsf{P}_{f\times f}\delta_{X}    \tag*{} \\
	&\ge \mathsf{P}_{\ang{\pi'_1,\pi'_2}}\delta_{X'}   \tag*{}  \\
	&:=\bfeqp{x'}{\pi}.  \tag*{}
\end{align}
We now prove condition \ref{item:weakelementarydoctrine-4}. Let $u:U\to A$ be a weak product $\weakbinprodg{r}{w}{u}{w}$ and let $h:U\to X\times X\times X\times X$ be the unique arrow induced by the projections $\pi_j\mathsf{r}_i:U\to X$, for $i,j=1,2$. If $t,t':u\to w$ are two arrows induced by the compositions $\pi_1\circ\mathsf{r}_1, \pi_2\circ\mathsf{r}_1$ and $\pi_1\circ\mathsf{r}_2, \pi_2\circ\mathsf{r}_2$, then we obtain condition \ref{item:weakelementarydoctrine-4} as follows
\begin{align}
	&(\mathsf{P}_{/A})_{\mathsf{r}_1}\bfeq{x}{\pi}\wedge (\mathsf{P}_{/A})_t\bfeq{x}{\pi}\wedge (\mathsf{P}_{/A})_{t'}\bfeq{x}{\pi} \tag*{} \\
	&:=\mathsf{P}_{h}(\mathsf{P}_{\ang{p_1,p_2}}\delta_{X}\wedge \mathsf{P}_{\ang{p_1,p_3}}\delta_{X}\wedge \mathsf{P}_{\ang{p_2,p_4}}\delta_{X}) \tag*{(\Cref{prop:equiv-elementarydoctrine}-\ref{item:equiv-elementarydoctrine-3})}\\
	&\le \mathsf{P}_h\mathsf{P}_{\ang{p_3,p_4}}\delta_{X} \tag*{}\\
	&=(\mathsf{P}_{/A})_{\mathsf{r}_2}\bfeq{x}{\pi}.\tag*{}
\end{align}
If $\Pdoc$ is a biased elementary doctrine, thanks to Lemma \ref{lemma:delta-is-determined-by-projections} we can repeat the above arguments to obtain that the slice doctrines of $\PP$ are biased. If $\PP$ has full weak comprehensions, the same holds for the slices of $\PP$.
\endexample

If we restrict to the 2-full subcategory $\textbf{BPD}_{\text{w.p.}}$ of \textbf{BPD} given by the biased primary doctrines whose base category has weak pullbacks, then we obtain the following characterization of slice doctrines.

\begin{lemma}\label{obs: slice doctrine as comma object}
		Let $\Pdoc$ be a biased elementary doctrine such that $\CC$ has weak pullbacks. For every object $A\in\CC$, the slice doctrine $\PP_{/A}$ is a strict comma object in the 2-category $\mathbf{BPD}_{\emph{\text{w.p.}}}$.
	\begin{proof}
	The argument follows the usual one for slice categories. Consider the terminal biased primary doctrine $\Pdocc{T}{\{\ast\}}$, which sends the terminal category to the singleton poset, and the 1-arrow $(A,\top_A):T\to \PP$, then the diagram
 \begin{equation}\label{diagram: strict commaobject}
 \begin{tikzcd}[column sep=large]
 	{\PP_{/A}} & \PP \\
 	T & \PP
 	\arrow["{1_\PP}", from=1-2, to=2-2]
 	\arrow["{(A,\top_A)}"', from=2-1, to=2-2]
 	\arrow["{(\ast,\top_\ast)}"', from=1-1, to=2-1]
 	\arrow["{(\mathsf{dom},id)}", from=1-1, to=1-2]
 	\arrow["\alpha", shorten <=11pt, shorten >=8pt, Rightarrow, from=1-2, to=2-1]
 \end{tikzcd} \end{equation}
where $\alpha$ is the 2-arrow defined $\alpha(f):=f$ for every $f:X\to A$ in $\CC/A$, has the following universal property. For every biased primary doctrine $\Pdocc{R}{\mathcal{D}}$ and pair of 1-arrows $(F,f):R\to P$ and $(G,g):R\to T$  in  $\textbf{BPD}_{\text{w.p.}}$ and 2-arrow $\beta: (F,f)\Rightarrow (A,\top_A)\circ(G,g)$, then there exists a unique 1-arrow $(H,h):R\to \PP_{/A}$ defined as $H(D):=\beta_D$ and $h_D:=f_D$  for every $D\in\mathcal{D}$, such that $(\mathsf{dom},id)\circ(H,h)= (F,f)$,  $(\ast,\top_\ast)\circ(H,h)= (G,g)$ and $\alpha(H,h)=\beta$.
\end{proof}\end{lemma}

\begin{rmk}\label{rmk: slice biased elementary and trict slices}\normalfont
Observe that in the proof of the above lemma the top horizontal arrow of diagram (\ref{diagram: strict commaobject}) does not preserve weak finite products. Hence, the result follows because of the choice of 1-arrows that we have made in \textbf{BPD}. By relaxing the notion of 1-arrows in \textbf{PD} to not require the preservation of finite products, a similar result can be obtained for the 2-full subcategory $\textbf{PD}_{\text{pb}}$ of \textbf{PD} given by the ordinary primary doctrines whose base category has pullbacks. Moreover, \Cref{obs: slice doctrine as comma object} restricts to (biased) elementary doctrines since the universal 1-arrow $(H,h)$ preserves (biased) fibered equalities in the sense of diagram (\ref{diagram: BED}) thanks to the fact that for every weak product $\wbinprod{D}{W}{D}$ the naturality of $\beta$ implies $\beta_DG(\pp_i)=\beta_W$ for $i=1,2$.
\end{rmk}

\begin{rmk}\label{obs:slice-doctrines}\normalfont
	 In many cases, elementary doctrines have weak equalizers (and weak pullbacks). For instance, if $\Pdoc$ has weak comprehensions and comprehensive diagonals, then $\CC$ has weak equalizers (and weak pullbacks), see \cite[Proposition 4.6]{maietti2013quotient}. In these cases, the slice categories have only weak finite products and then the slice doctrines are not in the realm of the elementary doctrines. The main examples of these situations are given by the elementary doctrines $F^{ML}$ and $G^{\textbf{mTT}}$ of Examples \ref{examples: elementary doctrines}. Their slices are biased elementary doctrines.

	  Intuitively, the slice doctrines arising from type theory take into account types depending on dependent types of a fixed closed type $A$. Indeed, every arrow $f:C\to A$ between closed types is a weak comprehension of a type $B(x)$ depending on $A$ (namely $B(x):=\underset{y:C}{\Sigma}{\mathsf{Id}_A(f(y),x)}$), and the fiber on $C$ is isomorphic to the fiber on $\underset{x:A}{\Sigma}B(x)$, which is given by the (equivalence classes of) types in the context $x:A,y:B(x)$. Moreover, for the projection arrows $\pi:\underset{x:A}{\Sigma}B(x)\to A$ the biased fibered equality is given by the identity type on $B(x)$ i.e.\ by the element $x:A, y_1:B, y_2:B\vdash \mathsf{Id}_{B(x)}(y_1,y_2)$ in the fiber on $\underset{x:A}{\Sigma}B(x)\times B(x)$.
	  
	  Hence, slice doctrines arising from type theories provide an explicit presentation  of types in non-empty context and their properties within the framework of doctrines. For instance, the arrow type $x:A\vdash B(x)\to C(x)$ of two dependent types can be described through a suitable weakened notion of exponential for these biased elementary doctrines, as done in \cite{cioffophd}.
      Moreover, the quotient completion of these doctrines are related to the interpretation of dependent types into the quotient model in \cite{maietti2009minimalist}, as discussed in \Cref{example: slice quotient mtt}.
	  
\end{rmk}

\example{(Weak subobjects)}\label{example:wed wek sub} \normalfont
	For a category with weak finite limits $\CC$, the functor of weak subobjects $\Psi_{\CC}$ of Example \ref{examples: elementary doctrines} is a biased elementary doctrine. Indeed, if $X$ is an object of $\CC$ and  $\weakbinprodg{p}{X}{W}{X}$ is a weak product, the element $\bfeq{X}{p}$ is given by the equivalence class, in the poset reflection of $\CC/W$, of a weak equalizer of the projections $\pp_1,\pp_2:W\to X$. The verification of conditions in Definition \ref{dfn:weakelementarydoctrine} follow from trivial weak pullback diagrams computations. This doctrine has full weak comprehensions.

\endexample


\section{Proof-irrelevant elements}\label{section:proof-irrelevant elements}

In this section, we detect particular classes of elements in the fibers of a biased elementary doctrine $\Pdoc$ on the weak products of $\CC$ that will be crucial for the rest of the article.

Given a list of objects $X_1,\dots,X_n\in\CC$ and a weak product $W$ with projections $\mathsf{p}_i:W\to X_i$, for $\lel{1}{i}{n}$, we can consider a commutative diagram of the form
\begin{equation}\label{diagram: proof-irrelevant of W}
\begin{tikzcd}[row sep=tiny]
	&&&&& {X_1} \\
	&& W &&& \vdots \\
	&& {W_1} &&& {X_n} \\
	U && \vdots \\
	&& {W_n} &&& {X_1} \\
	&& W &&& \vdots \\
	&&&&& {X_n}
	\arrow["{\mathsf{r_1}}", from=4-1, to=2-3]
	\arrow["{t_1}"', from=4-1, to=3-3]
	\arrow["{t_n}", from=4-1, to=5-3]
	\arrow["{\mathsf{r}_2}"', from=4-1, to=6-3]
	\arrow["{\mathsf{p}_n}"'{pos=0.6}, from=6-3, to=7-6]
	\arrow["{\mathsf{p}_1}"{pos=0.6}, from=2-3, to=1-6]
	\arrow["{\mathsf{q}^1_1}"'{pos=0.75}, from=3-3, to=1-6]
	\arrow["{\mathsf{p}_n}"'{pos=0.6}, from=2-3, to=3-6,crossing over]
	\arrow["{\mathsf{q}^n_2}"{pos=0.75}, from=5-3, to=7-6]
	\arrow["{\mathsf{q}^n_1}"{pos=0.75}, from=5-3, to=3-6]
	\arrow["{\mathsf{q}^1_2}"'{pos=0.75}, from=3-3, to=5-6,crossing over]
	\arrow["{\mathsf{p}_1}"{pos=0.6}, from=6-3, to=5-6, crossing over]
\end{tikzcd}
\end{equation}
where $\weakbinprodg{r}{W}{U}{W}$ and $\weakbinprodga{q}{X_i}{W_i}{X_i}{i}$, for $\lel{1}{i}{n}$, are weak binary products. For such a diagram we can consider the sub-poset of $\PP(W)$ given by the descent datum \[{\mathpzc{Des}}{(\mathsf{P}_{t_1}\bfeqa{X_1}{q}{1}\wedge\dots\wedge  \mathsf{P}_{t_n}\bfeqa{X_n}{q}{n})} \]
which consists of the elements $\alpha\in\PP(W)$ such that $$\mathsf{P}_{\mathsf{r}_1}\alpha\wedge \mathsf{P}_{t_1}\bfeqa{X_1}{q}{1}\wedge\dots\wedge  \mathsf{P}_{t_n}\bfeqa{X_n}{q}{n}\le \mathsf{P}_{\mathsf{r}_2}\alpha.$$
Apparently, the definition of the above elements of $\PP(W)$ depends on the particular commutative diagram (\ref{diagram: proof-irrelevant of W}) considered. The first result that we prove is that, surprisingly, they are completely determined by the weak product $W$.

	\begin{lemma}\label{lemma:the-main-prop-i}
	Let $\Pdoc$ be a biased elementary doctrine and let $\mathsf{p}_i:W\to X_i$, for $\lel{1}{i}{n}$, be a weak product of the objects $X_1,\dots,X_n\in\CC$. If we consider a commutative diagram as in (\ref{diagram: proof-irrelevant of W}) and another commutative diagram
	
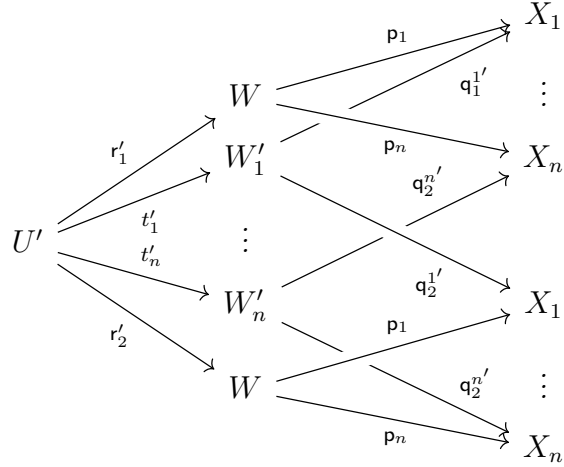
\begin{equation*}
		\begin{tikzcd}[row sep=tiny]
			&&&&& {X_1} \\
			&& W &&& \vdots \\
			&& {W'_1} &&& {X_n} \\
			U' && \vdots \\
			&& {W'_n} &&& {X_1} \\
			&& W &&& \vdots \\
			&&&&& {X_n}
			\arrow["{\mathsf{r}'_1}", from=4-1, to=2-3]
			\arrow["{t'_1}"', from=4-1, to=3-3]
			\arrow["{t'_n}", from=4-1, to=5-3]
			\arrow["{\mathsf{r}'_2}"', from=4-1, to=6-3]
			\arrow["{\mathsf{p}_n}"'{pos=0.6}, from=6-3, to=7-6]
			\arrow["{\mathsf{p}_1}"{pos=0.6}, from=2-3, to=1-6]
			\arrow["{\mathsf{q}^{1'}_1}"'{pos=0.75}, from=3-3, to=1-6]
			\arrow["{\mathsf{p}_n}"'{pos=0.6}, from=2-3, to=3-6,crossing over]
			\arrow["{\mathsf{q}^{n'}_2}"{pos=0.75}, from=5-3, to=7-6]
			\arrow["{\mathsf{q}^{n'}_2}"{pos=0.75}, from=5-3, to=3-6]
			\arrow["{\mathsf{q}^{1'}_2}"'{pos=0.75}, from=3-3, to=5-6,crossing over]
			\arrow["{\mathsf{p}_1}"{pos=0.6}, from=6-3, to=5-6, crossing over]
		\end{tikzcd}
	\end{equation*}
	where $\weakbinprodgp{r}{W}{U'}{W}$ and $\weakbinprodga{q}{X_i}{W'_i}{X_i}{i'}$, for $\lel{1}{i}{n}$, are  weak binary products, then we have the equality
	\[{\mathpzc{Des}}{(\mathsf{P}_{t_1}\bfeqa{X_1}{q}{1}\wedge\dots\wedge  \mathsf{P}_{t_n}\bfeqa{X_n}{q}{n})}={\mathpzc{Des}}{(\mathsf{P}_{t'_1}\bfeqa{X_1}{q}{1'}\wedge\dots\wedge  \mathsf{P}_{t'_n}\bfeqa{X_n}{q}{n'})}.\]
	\begin{proof}\normalfont
		The weak universal property of weak products induces arrows
		\begin{equation*}
			k:U'\to U \qquad h_{i}:W'_i\to W_i
		\end{equation*}
		such that $\mathsf{r}_jk=\mathsf{r}'_j$ and $\mathsf{q^i}_{j}h_{i}=\mathsf{q^i}'_{j}$ for $\lel{1}{i}{n}$ and $j=1,2$. We obtain the inclusion $``\subseteq"$ as follows. If $\mathsf{P}_{\mathsf{r}_1}\alpha\wedge \mathsf{P}_{t_1}(\bfeqa{X_1}{q}{1})\wedge\dots\wedge  \mathsf{P}_{t_n}(\bfeqa{X_n}{q}{n})\le \mathsf{P}_{\mathsf{r}_2}\alpha$, then 
		\begin{align}
			& \mathsf{P}_{\mathsf{r}'_1}\alpha\wedge \mathsf{P}_{t'_1}\bfeqa{X_1}{q}{1'}\wedge\dots\wedge  \mathsf{P}_{t'_n}\bfeqa{X_n}{q}{n'} \tag*{} \\
			& =\mathsf{P}_k\mathsf{P}_{\mathsf{r}_1}\alpha \wedge \mathsf{P}_{t'_1}\mathsf{P}_{h_{1}} \bfeqa{X_1}{q}{1}\wedge\dots\wedge  \mathsf{P}_{t'_n}\mathsf{P}_{h_{n}}\bfeqa{X_n}{q}{n}    \tag{\Cref{cor:delta-reindexings}} \\
			& = \mathsf{P}_k(\mathsf{P}_{\mathsf{r}_1}\alpha\wedge \mathsf{P}_{t_1}\bfeqa{X_1}{q}{1}\wedge\dots\wedge  \mathsf{P}_{t_n}\bfeqa{X_n}{q}{n})  \tag{\Cref{lemma:delta-is-determined-by-projections}}  \\
			& \le \mathsf{P}_k\mathsf{P}_{\mathsf{p}_2}\alpha \tag*{}\\
			& =\mathsf{P}_{\mathsf{r}'_2}\alpha \tag*{}.
		\end{align}
		The opposite inclusion $\supseteq$ follows similarly considering arrows $k':U\to U'$ and $h'_{i}:W_i\to W'_i$ such that $\mathsf{r}'_jk'=\mathsf{r}_j$ and $\mathsf{q^i}'_{j}h'_{i}=\mathsf{q^i}_{j}$ for $\lel{1}{i}{n}$, $j=1,2$.
	\end{proof}
\end{lemma}

Thanks to the above lemma, we can give the following definition.

\begin{dfn}\label{dfn:proof-irrelevant elements}\normalfont
Let $\Pdoc$ be a biased elementary doctrine. For every  weak product $\mathsf{p}_i:W\to X_i$, for $\lel{1}{i}{n}$, of the objects $X_1,\dots,X_n\in\CC$, we define the sub-poset of $\PP(W)$ of those elements $\alpha$ such that for some (and hence for all) commutative diagram as in (\ref{diagram: proof-irrelevant of W}) satisfies 
\begin{equation}\label{equation:proof-irrelevant}
	\mathsf{P}_{\mathsf{r}_1}\alpha\wedge \mathsf{P}_{t_1}\bfeqa{X_1}{q}{1}\wedge\dots\wedge  \mathsf{P}_{t_n}\bfeqa{X_n}{q}{n}\le \mathsf{P}_{\mathsf{r}_2}\alpha.
\end{equation}
We will refer to these elements as \textit{proof-irrelevant} elements (or \textit{strict predicates}) of the weak product $W$. 
\end{dfn}

Since, proof-irrelevant elements are closed under $\wedge$ and the top element $\top_W$ is trivially proof-irrelevant, we actually obtain a sub-inf-semilattice of $\PP(W)$. Observe that condition \ref{item:weakelementarydoctrine-4} of \Cref{dfn:weakelementarydoctrine} states that the proof-relevant equalities is a proof-irrelevant element. The adoption of the name proof-irrelevant will be clarified in \Cref{example: why-proof-irrelevant}.

	\begin{obs}\label{obs:proof-irrelevant of strict products}\normalfont
	If $\PP$ is a strict elementary doctrine and $X_1,\dots,X_n$ are objects of $\CC$, then (\ref{equation:delta of prod = prod of deltas}) implies that the proof-irrelevant elements of a strict product $X
	_1\times\dots\times X_n$ coincide with the whole fiber $\mathsf{P}(X
	_1\times\dots\times X_n)$.
\end{obs}

Proof-irrelevant elements are reindexed by projections.
\begin{prop}\label{prop:the-main-prop-ii}
	Let $\Pdoc$ be a biased elementary doctrine and let $X_1,\dots,X_n$ be objects of $\CC$. If $W$ is a weak product of $X_1,\dots,X_n$ with projections $\mathsf{p}_i:W\to X_i$, for $\lel{1}{i}{n}$, then the proof-irrelevant elements of $W$ are reindexed by projections.
	\begin{proof}
		Similar to the proof of Lemma \ref{lemma:delta-is-determined-by-projections}.
	\end{proof}
\end{prop}

\begin{rmk}\label{rmk:rbp=pi}\normalfont
	In general, we are not able to prove the converse of the above proposition, i.e.\ that every element $\alpha\in \PP(W)$ reindexed by projections is also proof-irrelevant according to (\ref{equation:proof-irrelevant}). However, an easy computation shows that in case $\PP$ has full weak comprehension and comprehensive diagonals, then the proof-irrelevant elements and the elements reindexed by projections of $W$ coincide.
\end{rmk}

In the following proposition we prove that two weak products of the same objects yield isomorphic proof-irrelevant elements.

	\begin{prop}\label{prop:the-main-prop-iii}
	Let $\Pdoc$ be a biased elementary doctrine and let $\mathsf{p}_i:W\to X_i$ and $\mathsf{p}'_i:W'\to X_i$, for $\lel{1}{i}{n}$, be two weak products of some objects $X_1,\dots,X_n\in\CC$. Then, the poset of proof-irrelevant elements of $W$ is isomorphic to the poset of proof-irrelevant elements of $W'$.
	\begin{proof}\normalfont
		We start fixing a commutative diagram as in (\ref{diagram: proof-irrelevant of W}). Now consider a similar diagram, for a weak product $\mathsf{p'}_i:W'\to X_i$, for $\lel{1}{i}{n}$, of the objects $X_1,\dots,X_n$, filled with 
	 weak products $\weakbinprodgp{r}{W'}{U'}{W'}$, $\weakbinprodga{q}{X_i}{W'_i}{X_i}{i'}$ and arrows $t'_i:U'\to W'_i$, for $\lel{1}{i}{n}$.
		Let
		\[{\mathpzc{Des}}{(\mathsf{P}_{t'_1}\bfeqa{X_1}{q}{1'}\wedge\dots\wedge  \mathsf{P}_{t'_n}\bfeqa{X_n}{q}{n'})}\subseteq \mathsf{P}(W')\] be the sub-poset of proof-irrelevant elements over $W'$.
		The weak universal property of weak products induces two arrows 
		\begin{equation*}
			h:W'\to W \qquad h':W\to W'
		\end{equation*}  
		such that $\mathsf{p}_ih=\mathsf{p'}_i$ and $\mathsf{p'}_ih'=\mathsf{p}_i$, for $\lel{1}{i}{n}$. We prove that the functors $\mathsf{P}_h$ and  $\mathsf{P}_{h'}$ provide an isomorphism between the sub-orders
		\begin{center}
			\begin{tikzcd}
				\mathsf{P}_h:{\mathpzc{Des}}{(\mathsf{P}_{t_1}\bfeqa{X_1}{q}{1}\wedge\dots\wedge  \mathsf{P}_{t_n}\bfeqa{X_n}{q}{n})}\arrow[shift left=2.1,""{name=U, below}]{r}{ } & {\mathpzc{Des}}{(\mathsf{P}_{t'_1}\bfeqa{X_1}{q}{1'}\wedge\dots\wedge  \mathsf{P}_{t'_n}\bfeqa{X_n}{q}{n'})}:\mathsf{P}_{h'}.\arrow[shift left=2.1,""{name=V, above}]{l}{ } \arrow[draw=none, from =U, to =V, "\cong" description]
			\end{tikzcd}
		\end{center}
		The weak universal property of weak products induces arrows
		\begin{equation*}
			k:U'\to U \qquad h_{i}:W'_i\to W_i
		\end{equation*}
		such that $\mathsf{r}_jk=h\mathsf{r}'_j$ and $\mathsf{q^i}_{j}h_{i}=\mathsf{q^i}'_{j}$ for $\lel{1}{i}{n}$, $j=1,2$. We first prove that if $\alpha$ is a proof-irrelevant element of $W$ then $\PP_h(\alpha)$ is a proof-irrelevant element of $\PP(W')$, in the following way: 
		\begin{align}
			&\mathsf{P}_{\mathsf{r}'_1}\mathsf{P}_h\alpha \wedge
			\mathsf{P}_{t'_1}\bfeqa{X_1}{q}{1'}\wedge\dots\wedge  \mathsf{P}_{t'_n}\bfeqa{X_n}{q}{n'} \tag*{} \\
			&= \mathsf{P}_k\mathsf{P}_{\mathsf{r}_1}\alpha \wedge \mathsf{P}_{t'_1}\mathsf{P}_{h_{1}}\bfeqa{X_1}{q}{1}\wedge\dots\wedge  \mathsf{P}_{t'_n}\mathsf{P}_{h_{n}}\bfeqa{X_n}{q}{n} \tag{\Cref{cor:delta-reindexings}} \\
			&=\mathsf{P}_k(\PP_{\mathsf{r}_1}\alpha \wedge\mathsf{P}_{t_1}\bfeqa{X_1}{q}{1}\wedge\dots\wedge  \mathsf{P}_{t_n}\bfeqa{X_n}{q}{n}) \tag{\Cref{lemma:delta-is-determined-by-projections}}\\
			&\le \mathsf{P}_k\mathsf{P}_{\mathsf{r}_2}\alpha \tag*{}\\
			&=\mathsf{P}_{\mathsf{r}'_2}\mathsf{P}_h\alpha. \tag*{}
		\end{align}
		A similar computation shows that $\mathsf{P}_{h'}$ sends proof-irrelevant elements of $W'$ into proof-irrelevant elements of $W$.
		The bijectivity follows from Proposition \ref{prop:the-main-prop-ii} since $\mathsf{p}_i(h\circ h')=\mathsf{p}_i\circ \text{id}$ and $\mathsf{p'}_i(h'\circ h)=\mathsf{p'}_i\circ \text{id}$ for $\lel{1}{i}{n}$.
	\end{proof}
\end{prop}

Thanks to Lemma \ref{lemma:the-main-prop-i}, Proposition \ref{prop:the-main-prop-ii}, and Proposition \ref{prop:the-main-prop-iii}, we can give the following definition.

\begin{dfn}\label{dfn:proof-irrelevant}\normalfont
	Let $\Pdoc$ be a biased elementary doctrine. For a list of objects $X_1,\dots,X_n\in\CC$, we can consider the poset obtained through the limit of the diagram of isomorphisms given by the restrictions of reindexing among the various presentations of proof irrelevant elements of the weak products of $X_1,\dots,X_n$ and denote it by $\Pstrict[X_1,\dots,X_n]$. We will refer to $\Pstrict[X_1,\dots,X_n]$ as \textit{proof-irrelevant elements} (or \textit{strict predicates}) of $X_1,\dots,X_n$.
\end{dfn}
In fact, an element in $\Pstrict[X_1,\dots,X_n]$ is the equivalence class of a proof-irrelevant element $\alpha\in\PP(W)$, of a weak product  $\mathsf{p}_i:W\to X_i$, for $\lel{1}{i}{n}$, up to the following equivalence relation: $\alpha$ is in relation with a proof irrelevant element $\alpha'\in\PP(W')$ of a weak product  $\mathsf{p}'_i:W'\to X_i$, for $\lel{1}{i}{n}$, if $\alpha=\PP_h(\alpha')$ for an arrow $h:W\to W'$, such that $\pp'_ih=\pp_i$, for $\lel{1}{i}{n}$.

We now describe proof-irrelevant elements for two examples of biased elementary doctrines. In particular, proof-irrelevant elements take their name from the following example.
	\begin{example}\label{example: why-proof-irrelevant}\normalfont
	Consider the elementary doctrine $\Pd{F^{ML}}{\ml}$  of and a closed type $A$. If $\lfloor f\rceil:X\to A$ is an object of $\ml/A$, then a canonical choice of weak product of $\lfloor f\rceil$ and $\lfloor f\rceil$ is given by the equivalence class of the common value of the two composites of the following weak pullback diagram
	\begin{center}
		\begin{tikzcd}
			W:=\underset{x_1,x_2:X}{\sum}\mathsf{Id}_A(f(x_1),f(x_2)) \arrow[]{r}{\pi_2} \arrow[swap]{d}{\pi_1}  & X \arrow[]{d}{f} \\
			X \arrow[swap]{r}{f} & A
		\end{tikzcd}
	\end{center}
	that we shall denote with $g:=f\pi_1:W\to A$. The slice doctrine ${F^{ML}}_{/A}$ sends the weak product $\lfloor w\rceil$ to  ${F^{ML}}_{/A}(\lfloor w\rceil):=F^{ML}(W)$, which is given by the equivalence classes of the dependent types on $W$:
	\[w:W\vdash B(w).\]
	 The fibered equality $\delta_{\lfloor f\rceil}$ is given by:
	\[w:W\vdash \mathsf{Id}_X(\pi_1(w),\pi_2(w)).\]
	Adopting the notation $x_i:= \pi_i(w):X$ for $i=1,2$ and $q:=\pi_3(w):\mathsf{Id}_A(f(x_1),f(x_2))$, the biased fibered equality $\delta_{\lfloor f\rceil}$ is simply given by the equality $\mathsf{Id}_X(x_1,x_2)$ of the first two components of $w$.

	Similarly, a weak product of $\lrfloor{g}$ and $\lrfloor{g}$ is given by the common value of the two composites of the following weak pullback
	\begin{center}
		\begin{tikzcd}
			U:=\underset{w,w':W}{\sum}\mathsf{Id}_A(g(w),g(w')) \arrow[]{r}{\pi_2} \arrow[swap]{d}{\pi_1}  & W \arrow[]{d}{g} \\
			W \arrow[swap]{r}{g} & A.
		\end{tikzcd}
	\end{center}
	that we shall denote with $h:=g\pi_1:U\to A$. The biased fibered equality $\delta_{\lfloor g\rceil}$ is an element of ${F^{ML}}_{/A}(h):=F^{ML}(U)$ given by:
	\[u:U\vdash \mathsf{Id}_W(\pi_1(u),\pi_2(u)).\]
	Instead, the element $\delta_{\lfloor f\rceil} \boxtimes \delta_{\lfloor f\rceil}$\ 
	corresponds to the dependent type
	\[u:U\vdash \mathsf{Id}_X(\pi_1(\pi_1(u)),\pi_2(\pi_1(u)))\times \mathsf{Id}_X(\pi_1(\pi_2(u)),\pi_2(\pi_2(u))).\]
	Adopting the notation $w:=\pi_1( u)$, $w':=\pi_2( u)$ and $p:=\pi_3(u):\mathsf{Id}_A(g(w),g(w'))$, then the element $\delta_{\lfloor f\rceil} \boxtimes \delta_{\lfloor f\rceil}$ corresponds to the component-wise equality
	
	\[\mathsf{Id}_X(x_1,x'_1)\times \mathsf{Id}_X(x_2,x'_2)\]
	
	which does not depend on the proof terms
	\begin{equation*}
		q: \mathsf{Id}_A(f(x_1),f(x_2)), \qquad
		q':\mathsf{Id}_A(f(x_2),f(x'_2)).
	\end{equation*} 
	Hence, we will refer to $\delta_{\lfloor f\rceil} \boxtimes \delta_{\lfloor f\rceil}$ as the\textit{ proof-irrelevant equality} of $W$, while we will refer to $\delta_{\lrfloor{g}}$ as  \textit{proof-relevant equality} of $W$.
	
	Similarly, we can describe proof-irrelevant elements of a weak product of some objects $f_i:X_i\to A$, for $i=1,\dots,n$, of $\ml/A$. A canonical choice of weak product of $f_1,\dots,f_n$ is given by the equivalence class of  $g:=\pi_1f_1:W\to A$ where \begin{equation}\label{equation:canonical-W}
		W:=\underset{x_1:X_1,\dots, x_n:X_n}{\sum}\mathsf{Id}_A(f_1(x_1),f_2(x_2))\times\dots\times\mathsf{Id}_A(f_{n-1}(x_{n-1}),f_n(x_n)).
	\end{equation} The proof-irrelevant elements of $W$ are types $w:W\vdash B(w)$ such that, if the type  $$\mathsf{Id}_X(x_1,x'_1)\times\dots\times \mathsf{Id}_X(x_n,x'_n)$$
	is inhabited, then $B(w)$ is inhabited if and only if $B(w')$ is inhabited. 	
\end{example}

\begin{example}\label{example:proof-irrelevant-weaksubobjects}\normalfont
	Let $\CC$ be a category with weak finite limits and consider the biased elementary doctrine $\Psi_{\CC}$ of weak subobjects of $\CC$. If $W$ is a weak product of a list of objects $X_1,\dots,X_n\in\CC$ with projections $\mathsf{p}_i:W\to X_i$, for $\lel{1}{i}{n}$, the unfolding of the descent condition (\ref{equation:proof-irrelevant}) of proof-irrelevant elements of $W$ involves several weak pullback diagrams computations. However, we can provide a concise description as follows. Let $\CC/(X_1,\dots,X_n)$ be the category of cones over $X_1,\dots,X_n$ and let ${\CC/(X_1,\dots,X_n)}\po$ be its poset reflection. We consider the functor 
	\begin{equation*}
		\M:\CC/(X_1,\dots,X_n)\po\to \Psi_{\CC}(W)
	\end{equation*}
	which sends the equivalence class of a cone $r_i:R\to X_i$, for $\lel{1}{i}{n}$ to the weak subobject given by the right dashed arrow of a weak limit of the following solid diagram
	
	\[\begin{tikzcd}
		& {R'} \\
		R && W \\
		{X_1} & \dotsm & {X_n}.
		\arrow[dashed, from=1-2, to=2-1]
		\arrow["\rho", dashed, from=1-2, to=2-3]
		\arrow["{r_1}"', from=2-1, to=3-1]
		\arrow["{r_n}"'{pos=0.2}, from=2-1, to=3-3]
		\arrow["{\mathsf{p}_n}", from=2-3, to=3-3]
		\arrow["{\mathsf{p}_1}"{pos=0.2}, from=2-3, to=3-1]
	\end{tikzcd}\]

In the opposite direction we have the functor
\begin{equation*}
	\mathcal{U}:\Psi_{\CC}(W)\to \CC/(X_1,\dots,X_n)\po
\end{equation*}
which sends a weak subobject $f:A\to W$ to the equivalence class of the cone given by the post compositions $\pp_i\circ f:A\to X_i$, for $\lel{1}{i}{n}$. In \cite[Proposition 4.5.2 and Theorem 4.5.3]{cioffophd}, we proved that the restriction of the functor $\mathcal{U}$ on the proof-irrelevant elements of $W$ is an isomorphism with inverse $\M$.

In particular, we observe that the biased fibered equality of an object $X\in\CC$ defined in \Cref{example:wed wek sub}, can be equivalently defined as the image of $\M$ on the trivial cone given by the two identities $1_X$, i.e. as the equivalence class of the right dashed arrow of a weak limit of the following solid diagram:
\begin{center}
	\begin{tikzcd}
		& D \arrow[dashed]{dl} \arrow[dashed]{dr}{\bfeq{X}{p}}&\\
		X \arrow[swap]{d}{1_X} \arrow[swap]{drr}[near start]{1_X} && W\arrow[]{dll}[near start]{\mathsf{p}_1} \arrow[]{d}{\mathsf{p}_2} \\
		X && X.
	\end{tikzcd}
\end{center}
\end{example}		
The following is an explicit description of proof-irrelevant elements for the weak subobjects doctrine of the category of free algebras $\setsG$.
\begin{example}\label{example: proof-irrelevant g-sets}\normalfont
Consider the category of free algebras $\setsG$ introduced in \Cref{examples: w prod } and the weak subobjects doctrine $\Pd{\Psi_{\setsG}}{\setsG}$. For a set X, the equalizer of the two projections of the weak products $X\times G\times X\times G$ is given by the arrow
\[\begin{tikzcd}
	E && {X\times G\times X\times G} && X
	\arrow["{\langle p_1,p_2\rangle}", shift left=2, squiggly, from=1-3, to=1-5]
	\arrow["{\langle p_3,p_4\rangle}"', shift right=2, squiggly, from=1-3, to=1-5]
	\arrow["{\langle m ,e\rangle}", squiggly, from=1-1, to=1-3]
\end{tikzcd}\]
where $(E,m)$ is the equalizer in $\set$ of $\langle p_1,p_2\rangle,\langle p_3,p_4\rangle$ and $e\in G$ is the neutral element of the group $G$. Hence, the biased fibered equality of $X$ is the weak subobject induced by the 5-tuples of the form $(x,g,x,g,k)\in X\times G\times X\times G\times G$, for $x\in X$ and $g,k\in G$. Now, given two sets $X,Y$ and their weak product $X\times G\times Y\times G$, its proof-relevant equality is the weak subobject induced by the tuples $(p,h,p',h',k)$ where $p,p'\in X\times G\times Y\times G$ and $h,h',k\in G$, such that 
\[\left\{ \begin{array}{ll}
p=p'\\
h=h'
\end{array}\right.\]
Instead, the proof-irrelevant equality is the weak subobject induced by the tuples $(p,h,p',h',k)$ such that
\begin{equation}\label{equation: proof irrelevant free algebras}\left\{ \begin{array}{lll}
	p_1=p_1' && p_3=p_3'\\
	p_2\cdot h=p_2'\cdot h'&& 	p_4\cdot h=p_4'\cdot h'
\end{array}\right.\end{equation}
Intuitively, an element $p$ of the weak product $X\times G\times Y\times G$ is given by a tuple $(p_1,p_2,p_3,p_4,h)$ and its $X,Y$ components are respectively given by $(p_1,p_2\cdot h)$ and $(p_3,p_4\cdot h)$. Conditions in (\ref{equation: proof irrelevant free algebras}) express the equality of the $X,Y$ components. Hence, a proof-irrelevant element of $\Psi_{\setsG}(X\times G\times Y\times G)$ is induced by a set $S$ of tuples closed by proof-irrelevant equality, i.e.\ if $(p_1,p_2,p_3,p_4,h)\in S$ and (\ref{equation: proof irrelevant free algebras}) holds, then $(p_1',p_2',p_3',p_4',h')\in S$.
\end{example}

 Implicational biased doctrines are defined as in Definition \ref{dfn:implicational}. A natural notion of existential and universal  biased doctrine is obtained requiring left and right adjoints as in Definitions \ref{dfn:existential-ED} and \ref{dfn:universal-ED}. However, we consider a weaker Beck-Chevalley condition and Frobenius reciprocity restricting to proof-irrelevant elements of a biased doctrine $\Pdoc$ as follows:
\begin{itemize}
	\item  the \textit{weak Beck-Chevalley} condition:
	for every arrow $f:Y\to X$ and any commutative diagram of the form
	\[\begin{tikzcd}[ column sep=small]
		& {X_1} && X_1 \\
		{V} && W \\
		& {Y} && X_2
		\arrow["{\mathsf{p}_1}", from=2-1, to=1-2]
		\arrow["{\mathsf{p}_2}"', from=2-1, to=3-2]
		\arrow["{\mathsf{p}_1}", from=2-3, to=1-4]
		\arrow["{\mathsf{p}_2}"', from=2-3, to=3-4]
		\arrow["1_{X_1}", from=1-2, to=1-4]
		\arrow["f"', from=3-2, to=3-4]
		\arrow["f'", from=2-1, to=2-3]
	\end{tikzcd}\]
	where $\wbinprod{X_1}{V}{Y}$ and $\wbinprod{X_1}{W}{X_2}$ are weak products, the canonical arrow $\exists_{\mathsf{p}_2}\circ \mathsf{P}_{f'}\alpha \le \mathsf{P}_{f}\circ \exists_{\mathsf{p}_2}\alpha$ is an isomorphism, for all proof-irrelevant $\alpha\in \mathsf{P}(W)$. Similarly, for the universal quantification.\\
	\item the \textit{weak Frobenius reciprocity}: for any projection $\mathsf{p}_i:W{\to} X_i$, element $\alpha\in \mathsf{P}(X_i)$, and $\beta\in \mathsf{P}(W)$ proof-irrelevant, the canonical arrow $\exists_{\mathsf{p}_i}(\mathsf{P}_{\mathsf{p}_i}\alpha \wedge \beta)\le\alpha\wedge\exists_{\mathsf{p}_i}\beta$ is an isomorphism.
\end{itemize}

\begin{rmk}\label{rmk: exist and univ for f}\normalfont
For existential elementary doctrines one can construct left adjoints along any arrow reindexing, see \cite[Remark 2.13]{maietti2013quotient}. The same holds for any existential biased elementary doctrine $\Pdoc$, i.e.\ for every arrow $f:Y\to X$ the functor $\mathsf{P}_f$ has left adjoint $\exists_f$ given by the functor which sends an element $\alpha\in \mathsf{P}(Y)$ to
\begin{equation}\label{equation:exists-f}
	\exists_f(\alpha):= \exists_{\mathsf{p}_2}(\mathsf{P}_{f'}\bfeq{p}{p}\wedge \mathsf{P}_{\mathsf{p}_1}\alpha),
\end{equation} 
where $f':K\to W$ is an arrow between weak products $\wbinprod{Y}{K}{X}$ and $\wbinprod{X}{W}{X}$, such that $\mathsf{p}_1\circ f'=f\circ \mathsf{p}_1$ and $\mathsf{p}_2\circ f'= \mathsf{p}_2$. 

In case of implicational and universal elementary doctrines one can construct also right adjoints along any arrow, see \cite[Remark 6.6]{maietti2013quotient}. The same holds for any implicational and universal biased elementary doctrine $\Pdoc$. The right adjoint along any arrow $f$ as above takes $\alpha\in \mathsf{P}(Y)$ to 
\begin{equation}\label{equation:univ-f}
	\forall_f(\alpha):=\forall_{\mathsf{p}_2}(\mathsf{P}_{f'}\bfeq{p}{p}\Rightarrow \mathsf{P}_{\mathsf{p}_1}(\alpha)).
\end{equation}
\end{rmk}

 \begin{example}\label{example: biased existential universal}\normalfont
	For a category with weak finite limits $\CC$ the doctrine $\Psi_{\CC}$ is existential. For every object $A\in\ml$, the doctrine ${F^{ML}}_{/A}$ is existential and universal. 
\end{example}

\section{Strictification}\label{section:strictification}
In this section, we construct an elementary doctrine from a biased one. In order to do that, we recall that the universal construction which freely adds strict products to a category $\CC$ is given by $(\mathsf{Fam}_{\mathsf{fin}}(\CC\op))\op$, see \cite{bunge1995symmetric}. The following is an explicit presentation of the product completion.

\begin{dfn}\label{dfn:product-completion}\normalfont
	Let $\CC$ be a category. The finite product completion of $\CC$ is the category $\CCs$ defined as follows:
	\begin{itemize}
		\item[-] {objects} of $\CCs$ are finite lists $[X_1,\dots,X_n]$ of objects of $\CC$.
		\item[-] arrows of $\CCs$ are pair $(f,s):[X_1,\dots,X_n]\to [Y_1,\dots,Y_m]$ such that $s:\{1,\dots,m\}\to \{1,\dots,n\}$ is a function and $f=[f_1,\dots,f_m]$ is a list of arrows $f_{i}:X_{s(i)}\to Y_{i}$ of $\CC$, for $\lel{1}{i}{m}$.
		The composition of two arrows $(f,s):[X_1,\dots,X_n]\to [Y_1,\dots,Y_m]$ and  $(g,t):[Y_1,\dots,Y_m]\to [Z_1,\dots,Z_k]$ is given by $$(g,t)\circ (f,s)= (g\ast f,s\circ t)$$
		where $g\ast f=[g_1\circ f_{t(1)},\dots,g_k\circ f_{t(k)}]$.
	\end{itemize} 
\end{dfn}
There is an obvious functor $S:\CC\to \CCs$, which sends an object $X\in\CC$ to the list $[X]\in\CCs$ and an arrow $f$ of $\CC$ to the arrow $(f, id_{\{1\}}):[X]\to [Y]$ of $\CCs$. The above construction provides a left bi-adjoint to the 2-forgetful functor $U:\textbf{Cart}\to \textbf{Cat}$, between the $2$-category \textbf{Cart} of small categories with strict finite products and functors preserving them and the 2-category of small categories \textbf{Cat}. The product of two lists is given by list concatenation. However, we observe that in case the category $\CC$ has weak finite products then the functor $S:\CC\to \CCs$ does not preserve them, neither it turns weak products into the strict ones.

	\begin{notation}\normalfont
 Given a category $\CC$, we will denote with $[f]:=S(f)$ the image through the functor $S:\CC \to \CCs$ of an arrow $f:X\to Y$ of $\CC$. Moreover, we will adopt the usual notation $\ang{[f],[g]}:[X]\to [Y,Z]$ to denote the unique arrow of $\CCs$ induced by the arrows $[f],[g]$ on the product $[Y,Z]$ of $\CCs$, i.e. $([f,g], c_1)$ where $c_1$ is the constant function equals to $1$. Finally, if $f_i:X_i\to Y_i$, for $\lel{1}{i}{n}$, are arrows of $\CC$, then we will denote with $[f_1,\dots,f_n]:[X_1,\dots,X_n]\to [Y_1,\dots,Y_n]$ the arrow  $([f_1,\dots,f_n],id_{\{1,\dots,n\}})$ of $\CCs$.
\end{notation}

We now prove that if $\Pdoc$ is a biased elementary doctrine then the  proof-irrelevant elements $\Pstrict([X_1,\dots,X_n])$  of a list of objects $X_1,\dots,X_n$ of $\CC$  induce a doctrine on the product completion of $\CC$. We start from the following construction.

\begin{prop}\label{prop:funcotiality-of-strictification}
	If $\Pdoc$ is a biased elementary doctrine, then we can define a functor 
	$\Pstrict:\CCs\op\to \mathsf{InfSL}$.
	\begin{proof}\normalfont
		The action of the  functor $\Pstrict$ on the non-empty lists in $\CCs$ is already defined. In particular, $\Pstrict([X])$ denotes nothing but the poset $\mathsf{P}(X)$. The terminal object is given by the empty list $[\ ]$ and the fiber $\Pstrict([\ ])$ is defined as the singleton set $\{\ast\}$. We impose that the reindexing along the unique (empty) arrow $!:[X_1,\dots,X_n]\to [\ ]$ sends $\ast$ to the top element of $\Pstrict([X_1,\dots,X_n])$.

		We now prove that the assignment is functorial.
		Let $(f,s):[X_1,\dots,X_n]\to [Y_1,\dots,Y_m]$ be an arrow of $\CCs$ and let $\pp_i:W\to X_i$, for $\lel{1}{i}{n}$ be a weak product of $X_1,\dots,X_n$ and $\mathsf{p}'_i:V\to Y_j$, for $\lel{1}{j}{m}$, be a weak product of $Y_1,\dots,Y_m$. Hence, we obtain an arrow $g:W\to V$ such that $\mathsf{p}'_i\circ g= f_i\circ\mathsf{p}_{s(i)}$, for $\lel{1}{i}{m}$. We now prove that if $\alpha$ is a proof-irrelevant element of $V$, then $\mathsf{P}_g(\alpha)$ is a proof-irrelevant element of $W$. In order to do that, we fix two commutative diagrams as in (\ref{diagram: proof-irrelevant of W}) for $W$ and $V$. We choose for $W$ some binary products $\weakbinprodg{r}{W}{U}{W}$ and $\weakbinprodg{q^i}{X_i}{W_i}{X_i}$ and arrows $t_i:U\to W_i$, for $\lel{1}{i}{n}$. We choose for $V$ some binary products $\weakbinprodgp{r}{V}{U'}{V}$ and $\weakbinprodgp{q^j}{Y_j}{V_j}{Y_j}$ and arrows $t'_j:U'\to V_j$, for $\lel{1}{j}{m}$. Now we consider commutative diagrams 
		of the form
		\[\begin{tikzcd}[row sep=tiny, column sep=small]
			& W && V &&& {X_{s(j)}} && {Y_j} \\
			U && {U'} &&& {W_{s(j)}} && {V_j} \\
			& W && V &&& {X_{s(j)}} && {Y_j}
			\arrow["{\mathsf{r}_1}", from=2-1, to=1-2]
			\arrow["{\mathsf{r}_2}"', from=2-1, to=3-2]
			\arrow["{\mathsf{r}'_1}"', from=2-3, to=1-4]
			\arrow["{\mathsf{r}'_2}", from=2-3, to=3-4]
			\arrow["g", from=1-2, to=1-4]
			\arrow["g"', from=3-2, to=3-4]
			\arrow["h", from=2-1, to=2-3]
			\arrow["{h_j}", from=2-6, to=2-8]
			\arrow["{\mathsf{q}_1^{s(j)}}", from=2-6, to=1-7]
			\arrow["{\mathsf{q}_2^{s(j)}}"', from=2-6, to=3-7]
			\arrow["{f_j}"', from=3-7, to=3-9]
			\arrow["{f_j}", from=1-7, to=1-9]
			\arrow["{\mathsf{q^j}'_1}"', from=2-8, to=1-9]
			\arrow["{\mathsf{q^j}'_2}", from=2-8, to=3-9]
		\end{tikzcd}\]
		for $\lel{1}{j}{m}$.
		Assuming $\mathsf{P}_{\mathsf{r}'_1}\alpha \wedge \mathsf{P}_{t'_1}\bfeqa{Y_1}{q}{1'}\wedge\dots\wedge  \mathsf{P}_{t'_n}\bfeqa{Y_n}{q}{n'}\le \mathsf{P}_{\mathsf{r}'_2}\alpha$ we obtain 
		\begin{align}
			&\mathsf{P}_{\mathsf{r}_1}\mathsf{P}_g\alpha\wedge \mathsf{P}_{t_1}\bfeqa{X_1}{q}{1}\wedge\dots\wedge  \mathsf{P}_{t_n}\bfeqa{X_n}{q}{n} \tag*{}\\
			&\le \PP_h\PP_{\mathsf{r}_1}\alpha\wedge \mathsf{P}_{t_{s(1)}}\mathsf{P}_{h_1}\bfeqa{Y_1}{q}{1'}\wedge\dots\wedge  \mathsf{P}_{t_{s(m)}}\mathsf{P}_{h_m}\bfeqa{Y_m}{q}{m'} \tag*{(\Cref{dfn:weakelementarydoctrine}-\ref{item:weakelementarydoctrine-3'})}\\
			&=\mathsf{P}_h(\mathsf{P}_{\mathsf{r}'_1}\alpha \wedge \mathsf{P}_{t'_1}\bfeqa{Y_1}{q}{1'}\wedge\dots\wedge  \mathsf{P}_{t'_n}\bfeqa{Y_n}{q}{n'}) \tag*{(\Cref{lemma:delta-is-determined-by-projections})}\\
			&\le \mathsf{P}_h\mathsf{P}_{\mathsf{r}'_2}\alpha \tag*{}\\
			&=\mathsf{P}_{\mathsf{r}_2}\mathsf{P}_g\alpha. \tag*{}
		\end{align}
		Hence, $\Pstrict(f,s)$ sends the isomorphism class of a proof-irrelevant element $\alpha$ of $X_1,\dots,X_n$ to the isomorphism class of the proof-irrelevant element $\mathsf{P}_g(\alpha)$ of $Y_1,\dots,Y_m$. The assignment of $\Pstrict$ is well defined and associative thanks to Proposition \ref{prop:the-main-prop-ii}.
	\end{proof}
\end{prop}
  \begin{thm}\label{thm:characterization-of-wed}
	If $\Pdoc$ is a biased elementary doctrine, then $\Pstrict$ is a strict elementary doctrine. Vice versa, for every strict elementary doctrine $\Pd{R}{\CCs}$, the pre-composition $\Pd{R\circ S}{\CC}$ is a biased elementary doctrine.
	\begin{proof}\normalfont
		Assuming that $\Pdoc$ is a biased elementary doctrine, we prove that $\Pd{\Pstrict}{\CCs}$ is a strict elementary doctrine as follows. For every object $[X_1,\dots,X_n]\in\CCs$ we consider a diagram as in (\ref{diagram: proof-irrelevant of W}) and define the fibered equality $\delta_{[X_1,\dots,X_n]}$ as the isomorphism class of the proof-irrelevant equality $\mathsf{P}_{t_1}\bfeqa{X_1}{q}{1}\wedge\dots\wedge  \mathsf{P}_{t_n}\bfeqa{X_n}{q}{n}$, which is in $\Pstrict[X_1,\dots,X_n,X_1,\dots,X_n]$. In particular, $\delta_{[X]}$ is the isomorphism class of the fibered equality $\bfeq{X}{p}$, of a weak product $\weakbinprodg{p}{X}{W}{X}$. Now we prove that $\delta_{[X_1,\dots,X_n]}$ satisfies conditions \ref{item:equiv-elementarydoctrine-1},\ref{item:equiv-elementarydoctrine-2} and \ref{item:equiv-elementarydoctrine-3} of Proposition \ref{prop:equiv-elementarydoctrine}.
		\ By Lemma \ref{lemma-conditions-3-4} and condition \ref{item:weakelementarydoctrine-1} of Definiton \ref{dfn:weakelementarydoctrine} we obtain $\top_{[X_1,\dots,X_n]}\le \Pstrict_{\Delta_{[X_1,\dots,X_n]}}\delta_{[X_1,\dots,X_n]}$. 
		Condition \ref{item:equiv-elementarydoctrine-2} is implicit in the definition of $\Pstrict$. Finally, condition \ref{item:equiv-elementarydoctrine-3} is obtained as follows. First consider a diagram as in (\ref{diagram: proof-irrelevant of W}) and a similar commutative diagram of weak finite products
		
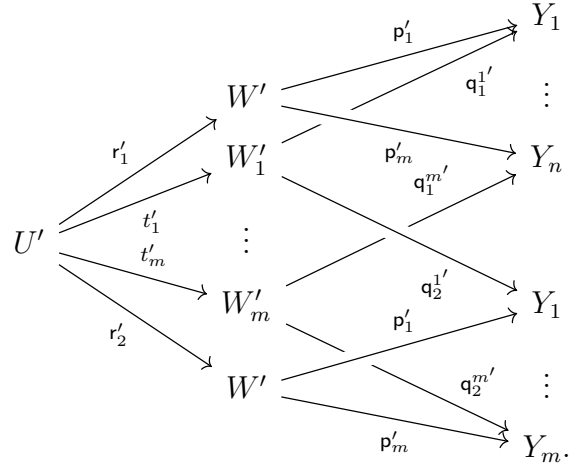
\begin{equation*}
			\begin{tikzcd}[row sep=tiny]
				&&&&& {Y_1} \\
				&& W' &&& \vdots \\
				&& {W'_1} &&& {Y_n} \\
				U' && \vdots \\
				&& {W'_m} &&& {Y_1} \\
				&& W' &&& \vdots \\
				&&&&& {Y_m.}
				\arrow["{\mathsf{r}'_1}", from=4-1, to=2-3]
				\arrow["{t'_1}"', from=4-1, to=3-3]
				\arrow["{t'_m}", from=4-1, to=5-3]
				\arrow["{\mathsf{r}'_2}"', from=4-1, to=6-3]
				\arrow["{\mathsf{p}'_m}"'{pos=0.6}, from=6-3, to=7-6]
				\arrow["{\mathsf{p}'_1}"{pos=0.6}, from=2-3, to=1-6]
				\arrow["{\mathsf{q}^{1'}_1}"'{pos=0.75}, from=3-3, to=1-6]
				\arrow["{\mathsf{p}'_m}"'{pos=0.6}, from=2-3, to=3-6,crossing over]
				\arrow["{\mathsf{q}^{m'}_2}"{pos=0.75}, from=5-3, to=7-6]
				\arrow["{\mathsf{q}^{m'}_1}"{pos=0.75}, from=5-3, to=3-6]
				\arrow["{\mathsf{q}^{1'}_2}"'{pos=0.75}, from=3-3, to=5-6,crossing over]
				\arrow["{\mathsf{p}'_1}"{pos=0.6}, from=6-3, to=5-6, crossing over]
			\end{tikzcd}
		\end{equation*}
		The element $\delta_{[X_1,\dots,X_n]}\boxtimes\delta_{[Y_1,\dots,Y_m]}$ is given by the isomorphism class of the element
		$$\mathsf{P}_{t}(\mathsf{P}_{t_1}\bfeqa{X_1}{p}{1}\wedge\dots\wedge  \mathsf{P}_{t_n}\bfeqa{X_n}{p}{n})\wedge \mathsf{P}_{t'}(\mathsf{P}_{t'_1}\bfeqa{Y_1}{p}{1'}\wedge\dots\wedge  \mathsf{P}_{{t'_m}}\bfeqa{Y_m}{p}{m'})$$
		for a commutative diagram 
		
		\[\begin{tikzcd}[row sep= tiny]
			&&&&& W \\
			&& V \\
			&& U &&& W' \\
			K \\
			&& U' &&& W \\
			&& V \\
			&&&&& W'
			\arrow["{\mathsf{s_1}}", from=4-1, to=2-3]
			\arrow["{t}"'{pos=0.65},from=4-1, to=3-3]
			\arrow["{t'}"{pos=0.65},from=4-1, to=5-3]
			\arrow["{\mathsf{s}_2}"', from=4-1, to=6-3]
			\arrow["{\mathsf{r}'_2}"{pos=0.75}, from=5-3, to=7-6]
			\arrow["{\mathsf{q}_1}"{pos=0.6}, from=6-3, to=5-6, crossing over]
			\arrow["{\mathsf{q}_2}"'{pos=0.6}, from=6-3, to=7-6]
			\arrow["{\mathsf{q}_1}"{pos=0.6}, from=2-3, to=1-6]
			\arrow["{\mathsf{r}_1}"'{pos=0.75}, from=3-3, to=1-6]
			\arrow["{\mathsf{q}_2}"'{pos=0.6}, from=2-3, to=3-6, crossing over]
			\arrow["{\mathsf{r}'_1}"{pos=0.75}, from=5-3, to=3-6]
			\arrow["{\mathsf{r}_2}"'{pos=0.75}, from=3-3, to=5-6, crossing over]
		\end{tikzcd}\]
		where $\weakbinprodg{q}{W}{V}{W'}$ and $\weakbinprodg{s}{V}{K}{V}$ are weak products.
		But, since $V$ is also a weak product of $X_1,\dots,X_n,Y_1,\dots,Y_m$, Lemma \ref{lemma:delta-is-determined-by-projections} implies the equality
		$$\delta_{[X_1,\dots,X_n]}\boxtimes\delta_{[Y_1,\dots,Y_m]}=\delta_{[X_1,\dots,X_n,Y_1,\dots,Y_m]}
		.$$\newline
		Now consider a strict elementary doctrine $\Pd{R}{\CCs}$ and the composition $\PP:=R\circ S$. We obtain that $\PP$ is a biased elementary doctrine setting for every weak product $\weakbinprodg{p}{X}{W}{X}$ the equality $\bfeq{X}{p}:=R_{\ang{[\mathsf{p}_1], [\mathsf{p}_2]}}\delta_{[X]}$.
	\end{proof}
\end{thm}

\begin{dfn}\label{dfn:strictification}\normalfont
	If $\Pdoc$ is a biased elementary doctrine, then the strict elementary doctrine $\Pd{\Pstrict}{\CCs}$ is called the \textit{strictification} of $\PP$.
\end{dfn}

\begin{rmk}\label{rmk: no universality strictification}\normalfont
The functor $S: \CC\to\CCs$, which does not preserve weak finite products, extends to a 1-arrow $(S,id):\PP\to \Pstrict$ between biased elementary doctrines. The strictification does not trivially provide a functor $\textbf{BED}\to\textbf{ED}$. This happens because if $(F,f):\PP\to R$ is a 1-arrow in \textbf{BED}, then we can trivially define a functor $F^\mathpzc{s}:\CCs \to \mathcal{D}_\mathpzc{s}$ which sends a list $[X_1,\dots,X_n]$ to the list $[F(X_1),\dots, F(X_n)]$. However, there is no obvious definition of the functors $f^\mathpzc{s}_{[X_1,\dots,X_n]}:\Pstrict([X_1,\dots,X_n])\to R^\mathpzc{s}([F(X_1),\dots, F(X_n)])$. A possible solution is to consider the subcategory of \textbf{BED} consisting of functors that preserve weak finite products, as a simple calculation shows that the functors $f_{(-)}$ preserve proof-irrelevant elements. However, even if the strictification provides a 2-functor under such a restriction,  the arrow $(S,id):\PP\to \Pstrict$ would still not provide the unit of a left biadjoint to the 2-forgetful functor {$U:\textbf{ED}\to \textbf{BED}$}. Finally, observe that if $\PP=R\circ S$ for an elementary doctrine $R\in\mathbf{ED}$, then we have an obvious 1-arrow $R\to \PP^\strict$. 
\end{rmk}

\begin{rmk}\label{rmk: rbp as sheaf condition}\normalfont
	As suggested by Rosolini, Definition \ref{dfn:reindex-determined-by-projections} is actually a sheaf condition of the topology on $\CCs$ generated by the (singleton) families of arrows of the form 
	\[\pp:[W]\to[X_1,\dots,X_n]\]
	where $p_i:W\to X_i$, for $\lel{1}{i}{n}$, is a weak product of the objects $X_1,\dots,X_n\in\CC$ and $\pp$ is the unique arrow induced by the projections $\pp_i$. Hence, in case of full comprehensions and comprehensive diagonals (see Remark \ref{rmk:rbp=pi}), the strictification of a biased elementary doctrine $\PP$ is actually a sheaf w.r.t.\ the above topology. The relationships between Grothendieck topologies, weak finite limits, and (biased) elementary doctrines are currently under investigation as part of ongoing work.
	
	\begin{obs}\label{obs: esistential proof irrelevant}\normalfont For a biased existential (universal) elementary doctrine,  we can consider three (or more) objects $X,Y,Z\in\CC$ and a weak product $\wbinprod{K}{U}{Z}$ of a weak product $\wbinprod{X}{K}{Y}$ and $Z$. Hence, we have a left (right) adjoint to the functor $\mathsf{P}_{\mathsf{p_2}}:\mathsf{P}(K)\to \mathsf{P}(U)$. Actually, the adjunctions restrict to proof-irrelevant elements providing a left (right) adjoint to the functor $\Pstrict_{p_2}$ as in the following diagram
	
	\[\begin{tikzcd}
		{\Pstrict[Z,X,Y]} & {\Pstrict[X,Y]} \\
		{\mathsf{P}(U)} & {\mathsf{P}(K).}
		\arrow[hook, from=1-1, to=2-1]
		\arrow[""{name=0, anchor=center, inner sep=0}, "{\exists_{\mathsf{p}_2}}", shift left=2, from=2-1, to=2-2]
		\arrow[""{name=1, anchor=center, inner sep=0}, "{\exists_{p_2}}", shift left=2, dashed, from=1-1, to=1-2]
		\arrow[hook, from=1-2, to=2-2]
		\arrow[""{name=2, anchor=center, inner sep=0}, "{\Pstrict_{p_2}}", shift left=2, from=1-2, to=1-1]
		\arrow[""{name=3, anchor=center, inner sep=0}, "{\mathsf{P}_{\mathsf{p_2}}}", shift left=2, from=2-2, to=2-1]
		\arrow["\dashv"{anchor=center, rotate=-90}, draw=none, from=1, to=2]
		\arrow["\dashv"{anchor=center, rotate=-90}, draw=none, from=0, to=3]
	\end{tikzcd}\]
Hence, we obtain the following property.
	\end{obs}
\end{rmk}
  \begin{prop}\label{dfn: existential universal biased}
	A biased elementary doctrine $\Pdoc$ is existential (implicational and universal) if and only if the strictification $\Pstrict$ is an existential (implicational and universal) elementary doctrine.
\begin{proof}\normalfont
If follows from Observation \ref{obs: esistential proof irrelevant} and the fact that existential doctrines have left adjoints to all reindexings, see \Cref{rmk: exist and univ for f}. Similarly, universal and implicational doctrines have right adjoint to all reindexings.
	\end{proof}
\end{prop}

\section{Quotient completion}\label{section: quotient completion}

	In this section, we present a quotient completion specifically designed for biased elementary doctrines. This construction mirrors the elementary quotient completion by Maietti and Rosolini \cite{maietti2013quotient}, extending its applicability to the more general context of biased elementary doctrines. To achieve this, proof-irrelevant elements will play a crucial role. We also demonstrate that this construction generalizes the exact completion of a category with weak finite limits \cite{carboni1998regular} in its more general form. Additionally, we discuss a universal property which exploits the notion of \textit{left covering} functor, introduced in \cite{carboni1998regular}, in the context of biased elementary doctrines. We mention that there is a different approach for generalizing the elementary quotient completion in \cite{dagnino2023quotients}, which may have non trivial intersection with our proposal.

We first note that the notion of $\PP$-equivalence relation, as defined in Definition \ref{dfn:P-eq relation}, relies on the properties of strict finite products. In case of a biased elementary doctrine, there is not a unique selection of arrows that define the properties of reflexivity, symmetry and transitivity for an element $\rho\in\PP(W)$, where $\weakbinprodg{p}{X}{W}{X}$ is a binary weak product. Therefore, in sight of Proposition \ref{prop:the-main-prop-ii}, a natural approach is to define $\PP$-equivalence relations as specific proof-irrelevant elements of $W$.

\begin{dfn}\label{dfn:P-eq relation for biased}\normalfont
	Let $\Pdoc$ be a biased elementary doctrine. A $\PP$-equivalence relation on an object $X\in\CC$ is an element $\rho\in\Pstrict[X,X]$ which is a $\Pstrict$-equivalence relation on $[X]\in\CCs$.
\end{dfn} 

	\begin{example}\label{example: PER as P-eq. rel }\normalfont
	In Example \ref{example:proof-irrelevant-weaksubobjects}, we have discussed proof-irrelevant elements of the biased elementary doctrine $\Psi_\CC$ of a category $\CC$ with weak finite limits. If $\weakbinprodg{p}{X}{W}{X}$ is a weak binary product, then $\Psi_{\CC}^\strict$ is given by the equivalence class of pair of arrows of the form
	\[r_1,r_2:R\to X.\]
	A straightforward verification demonstrates that the pair $r_1,r_2$ is a $\Psi_\CC$-equivalence relation on $X\in\CC$ if and only if it constitutes a pseudo equivalence relation, as defined in \cite{carboni1998regular}. Therefore, the classical definition of a pseudo equivalence relation in categories with weak finite limits (where finite products are also weak), is an instance of Definition \ref{dfn:P-eq relation for biased}  and, consequently, of proof-irrelevant elements. 
\end{example}

\begin{dfn}\normalfont
	Let $\Pdoc$ be a biased elementary doctrine and let $\rho$ be a $\PP$-equivalence relation on $X$.
	A \textit{quotient} of $\rho$ is an arrow $q:X\to C$ in $\CC$ such that $\rho\le \Pstrict_{[q]\times[q]}(\delta_{[C]})$ and, for every arrow $g:X\to Z$ such that $\rho\le \Pstrict_{[g]\times[g]}(\delta_{[Z]})$, there exists a unique arrow $h:C\to Z$ such that $g=h\circ q$. Stable, effective and effective descent quotients are defined as for elementary doctrines.
\end{dfn}

\begin{example}\label{example:quotients-wed-weaksubobject}\normalfont
	Let $\CC$ be a category with finite weak limits and let $X$ be an object of $\CC$. Example \ref{example: PER as P-eq. rel } implies that an arrow $q:X\to C$ is the coequalizer of a pseudo-equivalence relation $r_1,r_2:R\to X$ if and only if it is the quotient of the corresponding $\Psi_{\CC}$-equivalence relation $\lfloor\rho\rceil\in \Psi_{\CC}^{\strict}[X,X]$.
\end{example}

We now introduce the \textit{biased} elementary quotient completion, which parallels the construction by Maietti and Rosolini outlined in \cite{maietti2013quotient}. Let $\Pdoc$ be a biased elementary doctrine; we consider the category $\overline{\CC}$ whose 
\begin{itemize}
	\item[-]objects are pairs $(X,\rho)$ where $X$ is an object of $\CC$ and $\rho$ is a $\PP$-equivalence relation on $X$,
	\item[-] arrows $\lfloor f \rceil:(X,\rho)\to (Y,\sigma)$ are equivalence classes of arrows $f:X\to Y$ such that $\rho\le\nolinebreak \Pstrict_{[f]\times[f]}(\sigma)$. Two arrows $f,f'$ are equivalent when $\rho\le \Pstrict_{[f]\times[f']}(\sigma)$.
\end{itemize}
The assignment $\overline{\mathsf{P}}$ which sends an object $(X,\rho)\in\overline{\CC}$ to  $\overline{\mathsf{P}}(X,\rho):={\mathpzc{Des}}_{\rho}$ and an arrow $\lfloor f\rceil:(X,\rho)\to (Y,\sigma)$ to $\overline{\mathsf{P}}_{\lfloor f\rceil}:=\mathsf{P}_f$ is a well defined functor as it is stated in the following simple lemma. 
\begin{lemma}\label{lemma:weak-elementary-quotien-completion}
	If $\Pdoc$ is a biased elementary doctrine, then:
	\begin{enumerate}
		\item\label{item:weak-elementary-quotien-completion-ii} If $(X,\rho)$  and $(Y,\sigma)$  are two objects of $\overline{\CC}$ and $f:X\to Y$ is an arrow of $\CC$ such that $\rho\le \Pstrict_{[f]\times[f]}\sigma$, then $\mathsf{P}_f$ restricts to a map 
		\[\mathsf{P}_{f}: \mathpzc{Des}_{\sigma}\to \mathpzc{Des}_{\rho}.\]
		\item\label{item:weak-elementary-quotien-completion-iii} If $f,g:X\to Y$ are arrows of $\CC$ such that $\rho\le \Pstrict_{[f]\times[g]}\sigma$ and $\beta\in\mathpzc{Des}_{\sigma}$, then
		\[\mathsf{P}_f(\beta)=\mathsf{P}_g(\beta).\]
	\end{enumerate}
	\qed
\end{lemma}

We now prove that the quotient completion of a biased elementary doctrine leads to a strict elementary doctrine with suitable quotients. The following theorem appears as \cite[Lemma 5.3,5.4 and 5.5]{maietti2013quotient} for strict elementary doctrines.

\begin{thm}\label{thm:ed-of-weak-elementary-quotient-completion} 
	If $\Pdoc$ be a biased elementary doctrine, then $\Pdocc{\overline{\mathsf{P}}}{\overline{\CC}}$ is an elementary doctrine with effective quotients of effective descent. Moreover, if $\PP$ has weak full comprehensions and comprehensive diagonals, then $\overline{\PP}$ has strict full comprehensions, comprehensive diagonals and stable quotients.
	\begin{proof}\normalfont
		Given two objects $(X,\rho), (Y,\sigma)\in\overline{\CC}$, an easy verification shows that the strict products are given up to isomorphism  by
		\begin{equation*}
			(W,\rho\boxtimes\sigma),
		\end{equation*}
		where $W$ is a weak product of the objects $X,Y\in\CC$. The terminal object is given up to isomorphism by the pair $(1,\top_{[\ ]})$, where $1$ is a weak terminal object of $\CC$ and $\top_{[\ ]}$ is the unique element of $\Pstrict([\ ])$.
		Conditions \ref{item:equiv-elementarydoctrine-1}, \ref{item:equiv-elementarydoctrine-2} and \ref{item:equiv-elementarydoctrine-3} of Proposition \ref{prop:equiv-elementarydoctrine} are proved as in \Cref{thm:characterization-of-wed}. A $\overline{\mathsf{P}}$-equivalence relation $\mu$ on $(X,\rho)$ is a $\PP$-equivalence relation on $X$ such that $\rho\le\mu$. Hence, the quotient is given by $$\lfloor 1_X\rceil:(X,\rho)\to (X,\mu),$$
		and it is an effective quotient of effective descent. 
		
		Now assume that $\PP$ has full comprehensions and comprehensive diagonals. If $\alpha\in \overline{\mathsf{P}}(X,\rho)$, and $\llbrace\alpha\rrbrace:C\to X$ is a weak comprehension of $\alpha\in \mathsf{P}(X)$, then the strict full comprehension of $\alpha$ is given by 
		$$\lfloor \llbrace\alpha\rrbrace \rceil: (C,\rho')\to (X,\rho)$$ 
		where $\rho':=\Pstrict_{\lrfloor{\llbrace\alpha\rrbrace}\times\lrfloor{\llbrace\alpha\rrbrace}}\rho$. The diagonals are comprehensive by construction. We now prove that quotients are stable. Let $\lambda$ be a $\overline{\PP}$-eq. relation on the object $(Y,\sigma)$ and consider its quotient $\lrfloor{1_Y}:(Y,\sigma)\to(Y,\lambda)$. If $\lrfloor{f}:(X,\rho)\to(Y, \lambda)$ is an arrow, then the commutative diagram

		\[\begin{tikzcd}
			{(C, v)} && {(X,\rho)} \\
			& {(W,\sigma\boxtimes\rho)} \\
			{(Y,\sigma)} && {(Y,\lambda)}
			\arrow["{\lrfloor{\pi_2}}", from=1-1, to=1-3]
			\arrow["{\lrfloor{\pi_1}}"', from=1-1, to=3-1]
			\arrow["{\lrfloor{1_Y}}"', from=3-1, to=3-3]
			\arrow["{\lrfloor{c}}", from=1-1, to=2-2]
			\arrow["{\lrfloor{\mathsf{p}_1}}", from=2-2, to=3-1]
			\arrow["{\lrfloor{\mathsf{p}_2}}"', from=2-2, to=1-3]
			\arrow["{\lrfloor{f}}", from=1-3, to=3-3]
		\end{tikzcd}\]
		where $\wbinprod{Y}{W}{X}$ is a weak product, $c:={\lrbrace{\Pstrict_{\ang{\mathsf{p}_1,\mathsf{p}_2}}\Pstrict_{[1_Y]\times [f]}\lambda}}$ and $v:=\Pstrict_{[c]\times [c]}\Pstrict_{\ang{\mathsf{p}_1,\mathsf{p}_2}}\sigma\boxtimes\rho$ is a pullback diagram, thanks to the description of pullbacks through comprehensions of elementary doctrines with weak comprehensions and comprehensive diagonals (see \cite[Prop. 4.6]{maietti2013quotient}). We now prove that the element $(X,\rho)$ is isomorphic to the element $(C,w)$ where $w:=\Pstrict_{[c]\times [c]}\Pstrict_{\ang{\mathsf{p}_1,\mathsf{p}_2}}\lambda\boxtimes\rho$. If $h$ denotes an arrow of the form $X\to W$ induced by the arrows $f$ and $1_X$, such that $\mathsf{p_1}h=f$ and $\mathsf{p_2}h=1_X$, then since
		\[\top_X\le\Pstrict_{\ang{[f],[f]}}\lambda= \PP_h \Pstrict_{\ang{\mathsf{p}_1,\mathsf{p}_2}}\Pstrict_{[1_Y]\times [f]}\lambda \]
		there exists an arrow $g:X\to C$ such that $c\circ g=h$. Since
		\[\sigma\le \Pstrict_{[h]\times[h]}\Pstrict_{\ang{\mathsf{p}_1,\mathsf{p}_2}}\lambda\boxtimes\rho=\Pstrict_{[g]\times [g]}\Pstrict_{[c]\times [c]}\Pstrict_{\ang{\mathsf{p}_1,\mathsf{p}_2}}\lambda\boxtimes\rho=\Pstrict_{[g]\times [g]}w\]
		the arrow $g$ induces an arrow
		\[\lrfloor{g}:(X,\rho)\to (C,w).\]
		This arrow is the inverse of $\lrfloor{\pi_2}:(C,w)\to (X,\rho)$.
	\end{proof}
\end{thm}
\begin{obs}\normalfont
	Thanks to Observation \ref{obs:proof-irrelevant of strict products}, if $\Pdoc$ is an elementary doctrine, the elementary quotient completion and the biased one are isomorphic elementary doctrines.
\end{obs}
\begin{rmk}\label{rmk: exact completion that are not elementary quotient completion}\normalfont
	Examples \ref{example: PER as P-eq. rel } and \ref{example:quotients-wed-weaksubobject} imply that the exact completion of a category $\CC$ with weak finite limits is a particular case of the biased elementary quotient completion of $\Psi_\CC$, in the sense that there is an equivalence of the two categories $\overline{\CC}\cong \CC_{ex/wlex}$. In particular, the topos of $G$-sets $\Gsets$  and the exact completion of the slice categories of $\hotop$ introduced in  \Cref{examples: w prod } are instances of exact completions which are not elementary quotient completions of ordinary elementary doctrines due to the presence of weak finite products. We mention that in \cite{dagnino2023quotients} the authors provide a different approach to address this class of examples through what they call \textit{relational quotient completion}.
\end{rmk}

The following is our main example of biased elementary quotient completion which is neither an exact completion nor an elementary quotient completion of an ordinary elementary doctrine.

\begin{example}\label{example:main biased elementary quotient completion}\normalfont
	Given any closed collection $X\in\mathcal{CM}$, we can consider the slice doctrine $G^{\textbf{mTT}}_{/X}$ of the doctrine $G^{\textbf{mTT}}$ arising from the Minimalist Foundation, in Examples \ref{examples: elementary doctrines}. The biased elementary quotient completion $\overline{G^{\textbf{mTT}}_{/X}}$ gives rise to the category $\overline{\mathcal{CM}/X}$. This category plays a fundamental role in the quotient model of \cite{maietti2009minimalist}, as we will discuss in \Cref{example: slice quotient mtt} after providing a finer presentation of this category.

\end{example}

\begin{rmk}\label{rmk:equivalnet eqc}\normalfont
	We can add quotients to a biased elementary doctrine  $\Pdoc$ in a different way. Indeed, we could apply first the strictification and then the elementary quotient completion obtaining the elementary doctrine $\Pd{\overline{\PP^\strict}}{\overline{\CC^\strict}}$. However, this construction does not have the desired properties such as the universal property discussed. For instance, completing with quotients in this way, would imply that the projective objects of $\overline{\CC^\strict}$ would be more than the objects of $\CC$. An equivalent description of the biased elementary quotient completion in terms of the strictification is obtained observing that in the category $\overline{\CC^\strict}$ there is a collection of morphisms, that we will denote with $\Sigma$, made of the arrows of the form
	\begin{equation}\label{equation: Sigma arrows}
		\lrfloor{\pp}:([W], \PP^\strict_{\pp \times \pp}\sigma)\to ([X_1,\dots,X_n],\sigma)
	\end{equation}
	where $\pp_i:W\to X_i$, for $\lel{1}{i}{n}$, is a weak product and the arrow $\pp:[W]\to [X_1,\dots, X_n]$ is given by $\langle [\pp_1],\dots, [\pp_1]\rangle$. The class of morphisms $\Sigma$ form a \textit{calculus of right fractions} and, in order to obtain $\overline{\CC}$, we need to invert formally this class, see \cite{GabZis}, obtaining
	\[\overline{\CC}\cong \overline{\CC^\strict}[\Sigma^{-1}].\] Moreover, in case $\CC$ has a choice of weak finite products, $\overline{\CC}$ is a reflexive subcategory of $\overline{\CC^\strict}$. Indeed, there is an adjunction
	\[\begin{tikzcd}
		{[-]:\overline{\CC}} && {\overline{\CC^\strict}: L}
		\arrow[""{name=0, anchor=center, inner sep=0}, shift right=2, hook, from=1-1, to=1-3]
		\arrow[""{name=1, anchor=center, inner sep=0}, shift right=2, from=1-3, to=1-1]
		\arrow["\dashv"{anchor=center, rotate=-90}, draw=none, from=1, to=0]
	\end{tikzcd}\]
	where the functor $L([X_1,\dots, X_n], \rho):= (W, \PP^\strict_{\pp \times \pp}\rho)$, with $\pp_i:W\to X_i$, for $\lel{1}{i}{n}$, is the chosen weak product of $X_1,\dots, X_n\in\CC$. The action of $L$ on the arrows of $\overline{\CC^\strict}$ is defined similarly.
	
	Equivalently, the biased elementary quotient completion of $\PP$ can be obtained through the precomposition of $\overline{\PP^\strict}$ with the functor obtained by pulling-back $S:\CC\to \CCs$ along the forgetful functor $U:\overline{\CCs}\to \CCs$ which sends an object $([X_1,\dots, X_n], \rho)$ to its first component:
	\[\begin{tikzcd}
		{\overline{\CC}} & {\overline{\CCs}} & \mathsf{InfSL} \\
		\CC & \CCs
		\arrow["U", from=1-2, to=2-2]
		\arrow["{\overline{\PP^\strict}}", from=1-2, to=1-3]
		\arrow["S"', from=2-1, to=2-2]
		\arrow[from=1-1, to=2-1]
		\arrow["{S'}", from=1-1, to=1-2]
		\arrow["\lrcorner"{anchor=center, pos=0.125}, draw=none, from=1-1, to=2-2]
		\arrow["{\overline{\PP}}", bend left= 30pt, from=1-1, to=1-3]
	\end{tikzcd}\]
	
\end{rmk}

We now examine the universal property of the biased elementary quotient completion, which differs from the universal property stated in Theorem \ref{thm:elementary-quotient-completion-universal-property}. This distinction is not surprising, as a similar issue arises in the universal property of the exact completion for categories with weak finite limits. Indeed, as noted by Carboni and Vitale in \cite[\S 3.3]{carboni1998regular}, the exact completion construction for categories with weak finite limits yields a canonical functor"
\begin{equation*}
	\Gamma:\CC \to \CC_{ex}
\end{equation*}
that can not be the unit of a biadjunction between the 2-category of exact categories and exact functors, denoted by \textbf{EX}, and any definable 2-category of categories with weak finite limits, denoted by \textbf{WLEX}. However, the authors consider a special class of functors called \textit{left coverings}, namely those functors $F:\CC\to \mathcal{A}$ from a category $\CC$ with weak finite limits to an exact category $\mathcal{A}$ such that the image $F(L)$ of a  weak limit of a diagram $\mathcal{L}:D\to \CC$ factors through the limit of $F\mathcal{L}$ by a regular epimorphism, and provide a universal property of the exact completion in the sense of \cite[{Theorem 29}]{carboni1998regular} that we recall below. We refer to \cite{vitale1994left} for further details.

\begin{thm}{(Carboni and Vitale).}\label{thm:universality-ex/wlex}
	Let $\CC$ be a category with weak finite limits and let $\mathcal{A}$ be an exact category. The exact completion $\Gamma$ induces an equivalence between the category of left covering functors from $\CC$ to $\mathcal{A}$, and the category of exact functors from $\CC_{ex}$ to $\mathcal{A}$. The same holds for the regular completion, with respect to any regular category $\mathcal{A}$.\qed
\end{thm}

Taking advantages of the above result, we now define in the context of the biased elementary doctrines the analogous of left covering functors. From condition \ref{item:weakelementarydoctrine-3'} of Definition \ref{dfn:weakelementarydoctrine} we obtain a canonical 1-arrow 
$$(J,j):\mathsf{P}\to \overline{\mathsf{P}}$$
where the functor $J$ is defined on objects $X\in\CC$ as $J(X):=(X,\delta_{[X]})$ and on an arrow $f:X\to Y$ as $J(f):=\lfloor f\rceil:(X,\delta_{[X]})\to (Y,\delta_{[Y]})$. The functors $j_X$ are just the identities of $\mathsf{P}(X)$. Observe that $J$ does not preserve weak finite products and the pair $(J,j)$ defines a 1-arrow between biased elementary doctrines.

When $\PP$ has also weak comprehensions we can observe the following facts:
\begin{itemize}
	\item if $X,Y\in\CC$ and $\weakbinprodg{p}{X}{W}{X}$ is a weak product, then the unique arrow into the strict product
	$$
	\begin{tikzcd}
		(W,\delta_{[W]})\arrow{r} & (W, \delta_{[X]}\boxtimes \delta_{[Y]})
	\end{tikzcd}
	$$
	is a quotient of the $\overline{\mathsf{P}}$-equivalence relation $\delta_{[X]}\boxtimes \delta_{[Y]}$ over $(W,\delta_{[W]})$.
	\item If $\llbrace\alpha\rrbrace:X\to A$ is a weak comprehension of an element $\alpha\in \mathsf{P}(A)$, then $J(\llbrace\alpha\rrbrace)$ factors through
	the comprehension of $j(\alpha)$ via a $\overline{\mathsf{P}}$-quotient $$
	\begin{tikzcd}
		(X,\Pstrict_h\delta_{[A]})\arrow{r}{\llbrace \alpha\rrbrace} & (A,\delta_{[A]})\\
		(X,\delta_{[X]})\arrow[swap]{ur}{J(\llbrace\alpha\rrbrace)} \arrow[dashed]{u}
	\end{tikzcd}
	$$
	where $h:=[\llbrace\alpha\rrbrace]\times[\llbrace\alpha\rrbrace]:[X,X]\to [A,A]$.
	
\end{itemize}
The above observations lead to the following definition of left covering functors for biased elementary doctrine as follows.
\begin{dfn}\label{dfn:left-covering}\normalfont
	Let $\Pdoc$ be a biased elementary doctrine with weak full comprehensions and comprehensive diagonals and let $\Pd{R}{\mathcal{D}}$ be an object of \textbf{QD}. A pair $(F,f):\PP\to R$ is called \textit{left covering} when
	\begin{enumerate}
		\item\label{item:left-covering-1} The functor $F$ sends a weak product $W$ of the objects $X,Y\in\CC$ to the object $F(W)\in\mathcal{D}$ such that the unique arrow
		\begin{equation*}
			{\ang{F(\mathsf{p}_1),F(\mathsf{p}_2)}}: F(W)\to F(X)\times F(Y)
		\end{equation*}
		is a quotient of an $R$-equivalence relation.
		\item\label{item:left-covering-2} For every object $X\in\CC$, the functors $f_X:\PP(X)\to RF(X)$ preserve all the structure.  In particular, the functor $f_X$ preserves finite meets and for a weak product $\wbinprod{X}{W}{X}$ we have  $$f_{W}(\bfeq{W}{p})=R_{\ang{F(\mathsf{p}_1),F(\mathsf{p}_2)}}(\delta_{F(X)}).$$
		.
		\item\label{item:left-covering-3} if $\llbrace \alpha \rrbrace:X\to A$ is a weak comprehension of the element $\alpha\in \mathsf{P}(A)$, then the arrow $F(\llbrace \alpha\rrbrace):F(X) \to F(A)$ factors through $\llbrace f(\alpha) \rrbrace$ via a quotient of an $R$-equivalence relation.
	\end{enumerate}
\end{dfn}

\begin{rmk}\normalfont
	If $\rho\in \Pstrict[X,X]$ is a $\PP$-equivalence relation and we denote by the same $\rho$ a representative in $\PP(W)$, then $f_W\rho\in \mathpzc{Des}_{k_X}$, where $k_X$ is the $R$-kernel of $ \ang{F(\mathsf{p}_1),F(\mathsf{p}_2)}$, i.e.
	\begin{equation*}
		{R_{\ang{F(p_1),F(p_2)}\times \ang{F(p_1),F(p_2)}}\delta_{F(X)\times F(X)}}.
	\end{equation*}
Thus, by condition \ref{item:left-covering-1} of \Cref{dfn:left-covering} and the assumption that quotients are effective and of effective descent, we will simplify our notation by writing $f_{W}\rho\in R(F(X)\times F(X))$.
\end{rmk}
\begin{obs}\label{obs:lco-weak-subobjects}\normalfont
	Since (weak) equalizers can be built through (weak) comprehensions and comprehensive diagonals, \cite[Proposition 27]{carboni1998regular}   implies that the functor $F$ of a left covering 1-arrow $(F,f)$ is a left covering functor. Moreover,\cite[Lemma 21]{carboni1998regular} implies that $F$ preserves monomorphic families of arrows.
\end{obs}

We denote with $\textbf{Lco}(\PP,R)$ the category of left covering 1-arrows from $\PP$ to $R$ and natural transformations between them. In the following theorem we prove the universal property of the biased elementary quotient completion in style of Theorem \ref{thm:universality-ex/wlex}.

\begin{thm}\label{thm:univ-prop-weqc}
	Let $\PP$ be a biased elementary doctrine with full weak comprehensions and comprehensive diagonals, and let $R$ be an elementary doctrine in \emph{\textbf{QD}}. The pre-composition with the 1-arrow $(J,j)$ induces an equivalence between the following categories
	\begin{center}
		\begin{tikzcd}
			\emph{\textbf{QD}}(\overline{\mathsf{P}},R)\arrow[]{rrr}[name = UA, above]{(-)\circ(J,j)} &&& \emph{\textbf{Lco}}(\PP,R).
		\end{tikzcd}
	\end{center}
	\begin{proof}\normalfont
		We first prove that the functor $(-)\circ(J,j)$ is essentially surjective. Given a 1-arrow $(F,f)$ of $\textbf{Lco}(\PP,R)$ we can define a 1-arrow $(\bar{F},\bar{f}):\bar{\PP}\to R$ as follows. The functor $\bar{F}$ sends a projective object $(X,\delta_{[X]})$ to the image of $F$, i.e. $\bar{F}(X,\delta_{[X]}):=F(X)$. On the objects of the form $(X,\rho)$, the image $\bar{F}(X,\rho)$ is defined to be the codomain of the quotient of the $R$-equivalence relation $f_{W}\rho$, for a weak product $\weakbinprodg{p}{X}{W}{X}$. Similarly, if $\lfloor g\rceil:(X,\delta_{X})\to (Y,\delta_Y)$ is an arrow between projectives, then we define $\bar{F}(\lfloor g\rceil):=F(g)$. If  $\lfloor g\rceil:(X,\rho)\to (Y,\sigma)$ then we define $\bar{F}(\lfloor g\rceil)$ as the unique arrow induced by the quotients, which makes the following diagram commute
		\begin{center}
			\begin{tikzcd}
				\bar{F}(X,\delta_{X})\arrow[]{r}{F(g)} \arrow[]{d}&   \arrow[]{d}\bar{F}(Y,\delta_{Y}) \\
				\bar{F}(X,\rho)  \arrow[dashed]{r}{\bar{F}(\lrfloor{g})}  &   \bar{F}(Y,\sigma). 
			\end{tikzcd}
		\end{center}
		The functors $\bar{f}_{(-)}:\mathsf{P}(-)\to R(F(-))$ are defined as the functors $f$ on the projectives $(X,\delta_{[X]})$. On the elements $(X,\rho)$, it is a trivial verification to prove that the functor $f_X$ restricts to a functor
		\begin{equation*}
			\bar{f}_{(X,\rho)}:= f_X: \mathpzc{Des}{\rho}\to \mathpzc{Des}{f_{W}\rho}.
		\end{equation*}
		Since quotients are stable, $\bar{F}$ preserves finite products. We now prove that $(\bar{F},\bar{f})$ sends strict comprehensions to strict comprehensions. Indeed, as in the proof of Theorem \ref{thm:ed-of-weak-elementary-quotient-completion}, a strict comprehension of $\alpha\in \bar{P}(X,\rho)$ is given by 
		$$\lfloor \llbrace\alpha\rrbrace \rceil: (C,\rho')\to (X,\rho)$$ 
		where $\rho':=\Pstrict_{h}\rho$ and $h:=[\llbrace\alpha\rrbrace]\times[\llbrace\alpha\rrbrace]$. Since $R$ has strict comprehensions, it follows that a comprehension $\llbrace\bar{f}_{(X,\rho)}(\alpha)\rrbrace:D\to\bar{F}(X,\rho)$ of $\bar{f}_{(X,\rho)}(\alpha)$ is a monomorphism, hence $\delta_D= R_{\lrbrace{\bar{f}_{(X,\rho)}(\alpha)}\times\lrbrace{\bar{f}_{(X,\rho)}(\alpha)}}$ (see \cite[Corollary 4.8]{maietti2013quotient}). Condition \ref{item:left-covering-1} of Definition \ref{dfn:left-covering} implies that $D$ and $\bar{F}(C,\rho')$ are quotients of the same $R$-equivalence relation
		\begin{equation*}
			R_{F\llbrace\alpha\rrbrace\times F\llbrace\alpha\rrbrace}f_{W}\rho. 
		\end{equation*}
		Hence, we have proved that $(\bar{F},\bar{f})\in 	\textbf{QD}(\overline{\mathsf{P}},R)$. We now prove that the functor $(-)\circ(J,j)$ is fully faithful. Indeed if $(F,f)$ and $(G,g)$ are 1-arrows of $\textbf{Lco}(P,R)$ and $\theta:(F,f)\Rightarrow (G,g)$ is a 2-arrow, then it can be extended to a 2-arrow $\bar{\theta}:(\bar{F},\bar{f})\Rightarrow (\bar{G},\bar{g})$. The arrows $\bar{\theta}_{(X,\delta_{X})}$ on the projectives are defined as $\theta_X$. On the objects of the form $\bar{\theta}_{(X,\rho)}$ the arrow is defined as the unique arrow induced by quotients, which makes the following diagram commute
		\begin{center}
			\begin{tikzcd}
				\bar{F}(X,\delta_{X})  \arrow[]{d}\arrow[]{r}{\theta_X} &  \bar{G}(X,\delta_X) \arrow[]{d}\\
				\bar{F}(X,\rho)	 \arrow[dashed]{r}  &   \bar{G}(X,\rho). 
			\end{tikzcd}
		\end{center}
	\end{proof}
\end{thm}

\begin{rmk}\label{rmk: universal property biaes quotint completion}\normalfont
 As shown by a counterexample in \cite[\S 3.3]{carboni1998regular}, the composition of left covering functors is not necessarily left covering. Since left covering functors $F:\CC\to \mathcal{D}$ induce left covering 1-arrows between $\Psi_\CC$ and $\mathsf{Sub}_\mathcal{D}$, and also between $\mathsf{Sub}_\CC$ and $\mathsf{Sub}_\mathcal{D}$ if $\CC$ has finite limits (see \cite[Prop. 20]{carboni1998regular}), then repeating the same argument in \cite{carboni1998regular} it follows that also the composition of left covering 1-arrows is not necessarily left covering.
 Hence, \Cref{thm:elementary-quotient-completion-universal-property} does not provide a left biadjoint to the 2-forgetful functor $U:\textbf{QD }\to \textbf{BED}$.
\end{rmk}

	We conclude this section with some observations regarding the quotient completion of the slice doctrines. In particular, for a category with weak finite limits $\CC$, it is straightforward to prove that, for every object $A\in\CC$, the exact completion of the slice is equivalent to the slice of the exact completion in the following sense:
	\[(\CC/A)_{ex}\cong \CC_{ex}/(\Gamma A).\]
	We now demonstrate that, for biased elementary doctrines, the quotient completion and the slice construction commute when both are defined.

	 We first observe that if $\Pdoc$ is a biased elementary doctrine with weak comprehensions and comprehensive diagonals, then for every two arrows $x_i: X_i\to A$, $\lel{1}{i}{2}$, the following commutative diagram is a weak pullback 
	\begin{equation}\label{diagram: weak pullback}
		\begin{tikzcd}
	 	C && X_2 \\
	 	& W \\
	 	X_1 && Y,
	 	\arrow[dashed,"{\pi_2}", from=1-1, to=1-3]
	 	\arrow["x_2", from=1-3, to=3-3]
	 	\arrow[dashed, "{\pi_1}"', from=1-1, to=3-1]
	 	\arrow["x_1"', from=3-1, to=3-3]
	 	\arrow["\mathsf{p_2}"', from=2-2, to=1-3]
	 	\arrow["\mathsf{p_1}", from=2-2, to=3-1]
	 	\arrow["\llbrace\rho\rrbrace", from=1-1, to=2-2]
	 \end{tikzcd}\end{equation}
	 where $\weakbinprodg{\pp}{X_1}{W}{X_2}$ is a weak product and $\llbrace\rho\rrbrace$ is a weak comprehension of the element $\Pstrict_{[\pp]}\rho\in \mathsf{P}(W)$ with  $\rho:=\Pstrict_{[x_1]\times [x_2]}\delta_{[Y]}$. We will denote with $\lrbrace{\rho}_\strict$ the composition $$[\pp]\circ[\lrbrace{\rho}]:[C]\to[W]\to [X,X].$$
	 
	  The following lemma is a consequence of straightforward calculations; see also \cite[\S 3.6]{cioffophd}.
	 
	 \begin{lemma}\label{lemma:properties for slicing/quotientig}
Let $\Pdoc$ be a biased existential elementary doctrine with full weak comprehensions. For every $X\in\CC$ and $\rho\in\PP^\strict[X,X]$, the following hold:
\begin{itemize}
	\item $\exists_{\llbrace \alpha \rrbrace}\mathsf{P}_{\llbrace \alpha \rrbrace}\beta = \alpha \wedge \beta$, for every $\alpha,\beta\in\PP(X)$,
	\item  $\exists_{\fs{\llbrace \rho \rrbrace}}\Pstrict_{\fs{\llbrace \rho \rrbrace}}\gamma = \rho \wedge \gamma$, for every $\gamma\in\Pstrict[X,X]$.
	\item $\PP^\strict_{{\llbrace \rho \rrbrace}_\strict}\exists_{\fs{\llbrace \rho \rrbrace}}\beta=\beta$, for every $\beta\in\mathpzc{Des}(\PP^\strict_{{\llbrace \rho \rrbrace}_\strict \times {\llbrace \rho \rrbrace}_\strict}\delta_{[X,X]})$
\end{itemize}
\begin{proof}
The first equation follows from full comprehensions. The second one follows from the first and the fact that the counit of the adjunction $\exists_\pp\dashv \PP_\pp$ is an isomorphism. The third equation follows from \Cref{rmk: exist and univ for f} and Beck-Chevalley condition.
\end{proof}
\end{lemma}

	The following proposition relates the equivalence relations of the slice doctrines with the underlying doctrine.
	\begin{prop}\label{prop:relations in the slices}
Let $\Pdoc$ be a biased existential elementary doctrine with full weak comprehensions and comprehensive diagonals. If $x:X\to A$ is an arrow of $\CC$ and $\lrbrace{\rho}:C\to W$ is a weak comprehension of $\rho:=\Pstrict_{[x]\times[x]}\delta_{[A]}$, then
\begin{itemize}
	\item[(i)] if $\sigma$ is a $\PP$-equivalence relation on $X$ such that $\sigma\le \rho$, then $\Pstrict_{\lrbrace{\rho}_\strict}\sigma$ determines a $\mathsf{P}_{/A}$-equivalence relation on $x$ and $\mathpzc{Des}(\sigma)=\mathpzc{Des}({\Pstrict_{\lrbrace{\rho}_\strict}\sigma})$\footnote{$\mathpzc{Des}({\Pstrict_{\lrbrace{\rho}_\strict}\sigma})$ and also $\mathpzc{Des}({r})$ in $(ii)$ are taken w.r.t. the arrows $\pi_i:C\to X$ as in diagram (\ref{diagram: weak pullback}).},
	\item[(ii)] if $r\in\PP(C)$ determines a $\mathsf{P}_{/A}$-equivalence relation on $x$, then $\exists_{\lrbrace{\rho}_\strict}r$ is a $\PP$-equivalence relation on $X$ and $\mathpzc{Des}({r})=\mathpzc{Des}({\exists_{\lrbrace{\rho}_\strict}r})$.
\end{itemize}
\begin{proof}\normalfont
The reflexivity, symmetry and transitivity of $\Pstrict_{\fs{\llbrace \rho\rrbrace}}\sigma$ and $\exists_{\fs{\llbrace \rho \rrbrace}}r$ follow from that of $\sigma$ and $r$, through the Beck-Chevalley condition.
In $(i)$, it is trivial that  $\mathpzc{Des}(\sigma)\subseteq\mathpzc{Des}({\Pstrict_{\lrbrace{\rho}_\strict}\sigma})$. The opposite inclusion follows from Lemma \ref{lemma:properties for slicing/quotientig}, applying $\exists_{\lrbrace{\rho}_\strict}$ to $\Pstrict_{[\pi_1]}\alpha\wedge \Pstrict_{\fs{\llbrace \rho\rrbrace}}\sigma\le \Pstrict_{[\pi_1]}\alpha$.
In $(ii)$, it is trivial that $\mathpzc{Des}({\exists_{\lrbrace{\rho}_\strict}r})\subseteq\mathpzc{Des}({r})$. The opposite inclusion follows observing that $\PP^\strict_{{\llbrace \rho \rrbrace}_\strict}\exists_{\fs{\llbrace \rho \rrbrace}}r=r$ and $\exists_{\fs{\llbrace \rho \rrbrace}}r\le \rho$.
\end{proof}
	\end{prop}
	
	\begin{prop}\label{prop:slices-and-quotients-commute}
		If $\Pdoc$ is an existential biased elementary doctrine with weak full comprehensions and comprehensive diagonals, then for every object $A\in\CC$ 
		\begin{equation*}
			\overline{\mathsf{P}_{/A}}\cong \overline{\mathsf{P}}_{/(A,\delta_{[A]})}.
		\end{equation*}
	\begin{proof}
	Thanks to Lemma \ref{lemma:properties for slicing/quotientig} and Proposition \ref{prop:relations in the slices}, we can define a 1-arrow $(M,m):\overline{\mathsf{P}_{/A}}\to\overline{\mathsf{P}}_{/(A,\delta_{[A]})}$ of \textbf{QD} as follows. The functor $M$ sends an object $(x,r)$ of $\overline{\CC/A}$, with $x:X\to A$, to the object $$\lfloor x\rceil:(X,\exists_{\fs{\llbrace \rho\rrbrace}}r)\to (A,\delta_{[A]})$$ of $\overline{\CC}/(A,\delta_{[A]})$, where $\rho:=\Pstrict_{[x]\times[x]}\delta_{[A]}$. An arrow $f:(x,r)\to (y,s)$, where $y:Y\to A$, is sent to $\lfloor f\rceil:(X,\exists_{\fs{\llbrace \rho\rrbrace}}r)\to (Y,\exists_{\fs{\llbrace \sigma\rrbrace}}s)$, with $\sigma:=\Pstrict_{[y]\times[y]}\delta_{[A]}$. The natural transformation $m$ is given by the identity.
	
	We define an inverse 1-arrow $(N,n)$ to $(M,m)$ as follows. The functor $N$ sends an object $\lfloor x\rceil:(X,\lambda)\to (A,\delta_{[A]})$ of $\overline{\CC}/(A,\delta_{[A]})$ to the object 
	$$(x,\mathsf{P}_{\fs{\llbrace \rho\rrbrace}}\lambda)$$
	where $\rho:=\Pstrict_{[x]\times[x]}\delta_{[A]}$. An arrow $\lfloor f\rceil:\lfloor x\rceil\to \lfloor y\rceil$, with $\lfloor y\rceil:(Y,\sigma)\to (A,\delta_{[A]})$ is sent to the arrow $f$. The natural transformation $n$ is given by the identity.

	\end{proof}
	\end{prop}
	
	\begin{example}\label{example: slice quotient mtt}\normalfont
The primary examples of the above situation are the slice doctrines of $F^{ML}$ and $G^{\textbf{mTT}}$ as discussed in Examples \ref{examples: elementary doctrines}. Intuitively, in these cases, the isomorphism presented in \Cref{prop:slices-and-quotients-commute} indicates that the categories of families of setoids on a closed type $A$ in $\ml$, or in $\mathcal{CM}$ - which consist of pairs $(B(x), R)$, where $B(x)$ is a type depending on $A$ and $R$ is a dependent equivalence relation on $B(x)$ - are respectively equivalent to the slice categories $\overline{\ml}/(A,\mathsf{Id}_A)$ and $\overline{\mathcal{CM}}/(A,\mathsf{Id}_A)$. This fact is a particular instance of \cite[Proposition 4.12]{maietti2009minimalist}, which demonstrates that extensional dependent collections are related to the slice categories of $\overline{\mathcal{CM}}$. Furthermore, presenting the slice category of $\overline{\mathcal{CM}}/(A,\mathsf{Id}_A)$ as the elementary quotient completion of the doctrine $G^{\textbf{mTT}}_{/A}$ is beneficial for studying properties of $\overline{\mathcal{CM}}/(A,\mathsf{Id}_A)$, such as the cartesian closure. Indeed, as shown in \cite{cioffophd}, the existence of a weakened notion of exponentials in $G^{\textbf{mTT}}_{/A}$ implies the existence of exponentials in the base category $\overline{\mathcal{CM}}/(A,\mathsf{Id}_A)$ of the elementary quotient completion $\overline{G^{\textbf{mTT}}_{/A}}$. This provides an alternative approach for establishing the local cartesian closure of the quotient model presented in \cite{maietti2009minimalist}.
	\end{example}

\section{Concluding remarks}\label{section: concluding remarks}
 The properties of the exact completion of categories with weak finite limits have been extensively studied in the literature, as evidenced by several works, including \cite{menni2000exact, carboni2000locally, gran1998exact, emmenegger2020local}. These studies characterize categories for which the exact completion is extensive, (locally) cartesian closed, or a (quasi) topos. Maietti and Rosolini, and later Pasquali, generalized many of these results through the elementary quotient completion, as described in \cite{maietti2013quotient, mpr21}. In future work, we aim to present a comprehensive formulation of these properties using biased elementary doctrines.

For instance, as discussed in the previous sections, one of the main advantages of adopting biased elementary doctrines is their closure under the operation of slicing, see Example \ref{example:slicedoctrine-of-ed}. This property is particularly beneficial for generalizing theorems regarding the local cartesian closure of the exact completion to the context of the elementary quotient completion. Carboni and Rosolini showed in \cite{carboni2000locally} that the local cartesian closure of the exact completion of a category $\CC$ with weak finite limits is equivalent to assuming a weakened notion of exponential in $\CC$. However, Emmenegger pointed out in \cite{emmenegger2020local} that when $\CC$ has weak finite products, a different notion of exponential, called \textit{extensional} in \cite{emmenegger2020local}, may be assumed instead.

In the context of elementary doctrines, \cite[Proposition 6.7]{maietti2013quotient} presents the first result concerning the cartesian closure of the base category of the elementary quotient completion $\Pd{\overline{\PP}}{\overline{\CC}}$ for a suitable doctrine $\PP$. Additionally, two results regarding the local cartesian closure of $\overline{\CC}$ can be found in \cite{mpr21} and \cite{cioffophd}. In this context, the local cartesian closure is achieved through detailed computations in the slice categories, which may present some challenges. However, it is important to note that these results only extend those in \cite{carboni2000locally, emmenegger2020local} to the scenario where the base category$\CC$ has weak finite limits but strict finite products. To achieve a comprehensive generalization of these results that accounts for weak finite products, it is advantageous to work within the framework of biased elementary doctrines, as demonstrated in the second half of \cite[\S 3]{cioffophd}.

Under suitable conditions, if $\Pdoc$ is a universal and existential biased elementary doctrine, then the slice doctrines inherit these properties, becoming universal and existential as well \cite[Proposition 3.6.6]{cioffophd}. Furthermore, as proved in Proposition \ref{prop:slices-and-quotients-commute}, the operation of slicing and quotient commute.

Building on the work of \cite{emmenegger2020local}, we extended the notion of extensional exponentials for biased elementary doctrines and proved a generalization of \cite[Theorem 3.6]{emmenegger2020local} within this context, see \cite[Theorem 3.7.4]{cioffophd}. Consequently, the local cartesian closure of $\overline{\CC}$, is achieved by examining the cartesian closure of the quotients of the slice doctrines of $\CC$, see \cite[Theorem 3.7.5]{cioffophd}.

Finally, we noted in \Cref{rmk: no universality strictification} and \Cref{rmk: universal property biaes quotint completion} that neither the strictification nor the biased elementary quotient completion provides left biadjoints to the 2-forgetful functors \textbf{ED}$\to$\textbf{BED} and \textbf{QD}$\to$\textbf{BED}. However, we do not rule out the possibility of discovering a suitable 2-category for biased elementary doctrines that could yield such universal properties. We leave this investigation for future work, along with a potential analysis of 2-(co)monads for these constructions, similar to the approaches taken in \cite{Trotta_2020}, \cite{trotta2019existential} and \cite{emmpasqroscoalgebras} for various universal constructions within the context of elementary and primary doctrines.

\bibliography{bib2}
\end{document}